\documentclass[11pt]{article}
\usepackage[usenames]{color} %
\usepackage[tbtags]{amsmath}
\usepackage{amssymb,amsthm,amsfonts} %
\usepackage[utf8]{inputenc} %
\usepackage{hyperref, bm, bbm}
\usepackage{tikz}
\usepackage{comment}
\usepackage{listings}
\usepackage{booktabs}
\usepackage{caption}
\usepackage{colonequals}
\usepackage[paperwidth=8.5in,paperheight=11in,margin=1in]{geometry}
\usepackage{mathtools}
\usepackage{graphicx}
\usepackage{multicol}
\usepackage{multirow}
\usepackage{accents}
\usepackage{caption}
\usepackage{subfig}
\usepackage{mathrsfs, changepage}
\usepackage{enumitem,kantlipsum, lipsum}
\usepackage{tikz}
\usetikzlibrary{shapes.geometric}
\usepackage{siunitx}
\usetikzlibrary{calc}

\usepackage{enumitem}

\usepackage{fancyhdr}
\usepackage{lastpage}
\usetikzlibrary{arrows,arrows.meta, shapes.geometric}
\usetikzlibrary{shadows}
\usetikzlibrary{shadows.blur}

\usepackage{xcolor}

\DeclareRobustCommand\dout{\bgroup
 \markoverwith{\lower-0.6ex\hbox
 {\kern-0.02em\vbox{\hrule width.2em\kern0.35ex\hrule}\kern-.06em}}%
 \ULon}

\makeatletter
\newcommand{\distas}[1]{\mathbin{\overset{#1}{\kern\z@\sim}}}%
\newsavebox{\mybox}\newsavebox{\mysim}
\newcommand{\distras}[1]{%
  \savebox{\mybox}{\hbox{\kern3pt$\scriptstyle#1$\kern3pt}}%
  \savebox{\mysim}{\hbox{$\sim$}}%
  \mathbin{\overset{#1}{\kern\z@\resizebox{\wd\mybox}{\ht\mysim}{$\sim$}}}%
}

\makeatother
\makeatletter
\newtheoremstyle{problemstyle}  							%
        {3pt}                                               %
        {3pt}                                               %
        {\normalfont}                               		%
        {}                                                  %
        {\bfseries}                 						%
        {\normalfont\bfseries.}         					%
        {.5em}                                          	%
        {}

\theoremstyle{plain}
\newtheorem*{remark}{Remark}
\newcommand*{\rom}[1]{\expandafter\@slowromancap\romannumeral #1@}

\makeatletter
\newcommand{\mathleft}{\@fleqntrue\@mathmargin0pt}
\newcommand{\mathcenter}{\@fleqnfalse}
\makeatother

\newtheorem{theorem}{Theorem}[section]

\newtheorem{definition}[theorem]{Definition}

\newtheorem{proposition}[theorem]{Proposition}

\newtheorem{corollary}[theorem]{Corollary}

\newtheorem{lemma}[theorem]{Lemma}

\def\E{\mathbb{E}}
\def\Ex{\mathop{\mathbb{E}}}
\newcommand{\Prx}{\mathop{\mathsf{Pr}}}
\def\Cov{\mathrm{Cov}}
\def\Covx{\mathop{\mathrm{Cov}}}
\def\Var{\mathrm{Var}}
\def\Varx{\mathop{\mathrm{Var}}}

\def\R{\mathbb{R}}

\def\Q{\mathbb{Q}}
\def\P{\mathbb{P}}

\def\E{\mathbb{E}}

\def\Law{\text{Law}}

\usepackage[mathscr]{euscript}

\newcommand{\comp}{\textnormal{comp}}
\newcommand{\stat}{\textnormal{stat}}

\renewcommand{\epsilon}{\varepsilon}

\def\Part{\mathrm{Part}}

\definecolor{cblack}{rgb}{0,0,0}
\definecolor{cblue}{rgb}{0.121569,0.466667,0.705882}    %
\definecolor{corange}{rgb}{1.000000,0.498039,0.054902}  %
\definecolor{cgreen}{rgb}{0.172549,0.627451,0.172549}   %
\definecolor{cred}{rgb}{0.839216,0.152941,0.156863}     %
\definecolor{cpurple}{rgb}{0.580392,0.403922,0.741176}  %
\definecolor{cbrown}{rgb}{0.549020,0.337255,0.294118}   %
\definecolor{cpink}{rgb}{0.890196,0.466667,0.760784}
\definecolor{cgray}{rgb}{0.498039,0.498039,0.498039}
\definecolor{cgreen2}{rgb}{0.7372549019607844, 0.7411764705882353, 0.13333333333333333}

\hypersetup{
  linkcolor  = cblue,
  citecolor  = cgreen,
  urlcolor   = corange,
  colorlinks = true,
}

\title{Statistical-to-computational gaps in central limit theorems for high-dimensional tensors}
\date{November 11, 2026}

\usepackage{authblk}
\author{Debsurya De\thanks{Email: \texttt{dde4@jhu.edu}}\,\,}
\author{Dmitriy Kunisky\thanks{Email: \texttt{kunisky@jhu.edu}.}}
\affil{Department of Applied Mathematics \& Statistics, Johns Hopkins University}

\begin{document}

\maketitle

\thispagestyle{empty}

\begin{abstract}
    We study high-dimensional central limit theorems (CLTs) as hypothesis testing problems, comparing the number of summands needed for statistical and computational indistinguishability of the normalized sum from the limiting Gaussian distribution as the dimension and number of summands grow together.
    For the off-diagonal entries of sums of i.i.d.\ rank-one tensors formed from standard Gaussian vectors (equivalently, the problem of distinguishing Wigner from Wishart tensors), we prove total variation bounds sharpening previous estimates of Mikulincer~(2020) and give matching polynomial-time tests.
    This shows that the statistical and computational thresholds coincide for this CLT, contrary to a proposal by Kunisky, Moore, and Wein~(2024).
    On the other hand, we also show that, for tensors of order at least three, sufficient Gaussian dilution of the summands by adding independent Gaussian noise and rescaling gives a CLT model with a statistical-to-computational gap (in the above setting and ones with more general rank-one summand distributions). To the best of our knowledge, such a phenomenon has not been observed before in any high-dimensional CLT. In particular, in the relevant parameter regimes, we show that a brute-force search achieves strong detection, while polynomials of degree up to $n^{\Omega(1)}$ cannot achieve strong separation. The latter lower bounds control low-degree $\chi^2$-divergences by a novel method based on cumulant generating functions, and we further illustrate the use of these tools by applying them to masked tensor CLT models where only a prescribed subset of tensor entries is revealed.
\end{abstract}

\clearpage

\tableofcontents

\pagestyle{empty}

\clearpage

\setcounter{page}{1}
\pagestyle{plain}

\section{Introduction}

The multivariate central limit theorem (CLT) is a centerpiece of classical probability theory, stating that if $X_1, X_2, \dots \sim \mu$ are i.i.d.\ random vectors for $\mu$ a probability measure on $\R^N$ that is \emph{centered} (having mean zero) and \emph{isotropic} (having covariance the identity matrix $I_N$), then we have the convergence in distribution
\[ r^{-1/2} \sum_{i = 1}^r X_i \Rightarrow Z \]
for fixed $N$ as $r \to \infty$, for a standard Gaussian vector $Z \sim \mathcal{N}(0, I_N)$.
Here and below we abbreviate
\[ \mathcal{N} \colonequals \mathcal{N}(0, I_N), \]
leaving $N$ implicit.
We call $\mu$ in such settings the \emph{summand distribution}.
Equivalently, writing $*$ for the convolution of probability measures and $c \cdot \mu$ for the scaling of a probability measure $\mu$ by a deterministic scalar $c \in \R$ (the law of $cX$ when $X \sim \mu$), this is the weak convergence
\[ r^{-1/2} \cdot \mu^{*r} \xrightarrow{\mathrm{(w)}} \mathcal{N}. \]
A long line of work has also studied non-asymptotic versions of such statements, controlling distributional distances such as $\chi^2$-divergence, total variation, Wasserstein, and others, $d(r^{-1/2} \cdot \mu^{*r}, \mathcal{N})$, as a function of $n$ for various choices of natural sequences of measures $\mu = \mu_n$ over $\R^N$ for some growing dimension $N = N(n)$ and growing number of summands $r = r(n)$.
For broad treatments, see the monograph on Stein's method~\cite{CGS-2011-NormalApproximationStein} and the survey of high-dimensional CLTs over rectangles~\cite{CCKK-2023-HighDimensionalDataBootstrap}.
Specific results include bounds over convex sets~\cite{Bentkus-2003-BerryEsseenDimension}, Wasserstein bounds~\cite{Zhai-2018-HighDimensionalCLTW2,Bonis-2020-MultivariateCLTWasserstein}, and entropy bounds~\cite{EMZ-2018-HighDimensionalCLTMartingale}.
For the classical one-dimensional CLT, quantitative bounds in total variation and relative entropy are given in~\cite{BCG-2014-EntropicBerryEsseen}, while convergence in $\chi^2$-divergence and more general R\'enyi divergences is studied in~\cite{BCG-2019-RenyiDivergenceCLT}.

\subsection{Testing for CLTs}

We take a \emph{hypothesis testing} perspective on such questions: if for instance $d = d_{\mathrm{TV}}$ is the total variation (TV) distance, and we consider a sequence of $\mu = \mu_n$ over $\R^N$ and $r = r(n)$, then the question of whether $d_{\mathrm{TV}}(r^{-1/2} \cdot \mu^{*r}, \mathcal{N})$ converges to 1 as $n \to \infty$ is precisely the question of whether there is an arbitrary \emph{test} $t = t_n: \R^N \to \{0, 1\}$ that can with high probability \emph{strongly distinguish} between the distributions $r^{-1/2} \cdot \mu^{*r}$ and $\mathcal{N}$, i.e., one that can achieve
\[ \lim_{n \to \infty} \,\,\Prx_{Y \sim r^{-1/2} \cdot \mu^{*r}}\,\,[t_n(Y) = 1] = \lim_{n \to \infty} \,\,\Prx_{Y \sim \mathcal{N}}\,\,[t_n(Y) = 0] = 1. \]

Following recent literature in high-dimensional statistics~\cite{ZK-2016-Review,GMZ-2022-DisorderedSystemsHardness}, the aim of this paper is to study the difference between the \emph{statistical} complexity of high-dimensional CLTs as testing problems (where arbitrary tests $t_n$ are allowed above) and their \emph{computational} complexity (where only computationally efficient tests $t_n$ are allowed, for instance ones computable in polynomial time in $n$).
For the latter, we consider the average-case complexity of CLTs through the lens of \emph{low-degree polynomial} algorithms~\cite{kunisky2019notes,Wein-2025-LowDegreeSurvey}, as detailed in Section~\ref{sec:low-degree}.
Since as $r \to \infty$ the distribution $r^{-1/2} \cdot \mu^{*r}$ becomes closer to $\mathcal{N}$, we expect testing to become harder as $r = r(n)$ grows faster.
Thus, we phrase our results in terms of statistical and computational \emph{thresholds} for $r$, aiming to identify scalings $r_{\stat}(n)$ and $r_{\comp}(n)$ such that arbitrary tests can distinguish the two distributions in the above sense if $r \ll r_{\stat}$ and cannot if $r \gg r_{\stat}$, and polynomial-time tests can distinguish if $r \ll r_{\comp}$ and cannot if $r \gg r_{\comp}$.
(In the examples we study we will consider $N \asymp n^p$ for some $p \geq 1$ and $r(n) \sim n^{\alpha}$ for various $\alpha > 0$, so the above asymptotics will be differences on this polynomial scale in the ambient parameter $n$.)

\subsection{Statistical-to-computational gaps and the two CLTs}

In particular, our starting point is a tantalizing proposal in the recent work \cite{KMW-2024-TensorCumulantsInvariantInference} that certain high-dimensional CLTs might exhibit \emph{statistical-to-computational gaps}.
This is the phenomenon where the statistical and computational thresholds differ macroscopically, in our case having
\[ r_{\comp} \ll r_{\stat}. \]
We note that, if this happened, it would raise an interesting philosophical quandary: when a high-dimensional CLT has this behavior for a sequence of $\mu_n$, there are really \emph{two} different CLTs.
First, for a slower scaling of the number of summands $r(n) \gg r_{\comp}(n)$, the CLT ``kicks in'' for efficiently computable test statistics.
Then, only for the larger number of summands $r(n) \gg r_{\stat}(n)$ does the CLT kick in for all test statistics.
It seems perhaps application-specific which of these is more meaningful.
While previous work has almost exclusively focused on the \emph{statistical CLT} determined by $r_{\stat}(n)$, which is perhaps more natural mathematically, certainly for applications to algorithms the \emph{computational CLT} determined by $r_{\comp}(n)$ is more important, determining whether efficient computations can ``see'' that the CLT has not yet kicked in (similar for instance to the notion of \emph{computational indistinguishability} in cryptography \cite[Section~3.2]{Goldreich-2001-FoundationsCryptography}).
Further, if such a CLT had a role in a physical argument or calculation, it is an intriguing question which CLT would be more relevant, perhaps related to long-studied questions about what kinds of abstract computations can be physically realized in a short time~\cite{Moore-1998-FiniteDimensionalAnalogComputers,Moore-2002-ComputationalComplexityPhysics,Aaronson-2005-NPCompletePhysicalReality,Aaronson-2013-PhilosophersShouldCare,BP-2021-AnalogModelsComputation}.

However, to the best of our knowledge no prior work has actually established strong evidence that any high-dimensional CLT exhibits such a gap.
The work \cite{KMW-2024-TensorCumulantsInvariantInference} merely suggested that this might be the case, establishing (in a certain sense involving low-degree polynomial algorithms) the scaling of $r_{\comp}(n)$ for the model they study but relying only on a proposal that a statistical bound from \cite{Mikulincer-2020-CLTWishartTensors} might be sharp.
That work considers high-dimensional CLTs in the Wasserstein distance rather than the TV distance, in a different model than \cite{KMW-2024-TensorCumulantsInvariantInference}.
Its bounds give sufficient conditions for convergence to Gaussianity, without matching lower bounds establishing statistical distinguishability in the proposed gap regime.
As we detail below in Section~\ref{sec:tensor-clt}, \cite{Mikulincer-2020-CLTWishartTensors} considers a model where $\mu_n$ is the law of $X = (v_{i_1} \cdots v_{i_p})_{1 \leq i_1 < \cdots < i_p \leq n}$ for $v \in \R^n$ a random vector with i.i.d.\ entries (under suitable regularity assumptions), which may be interpreted as a problem of testing between ``Wigner'' and ``Wishart'' models of random tensors.
Because of technical limitations of their approach, \cite{KMW-2024-TensorCumulantsInvariantInference} only study a certain proxy for this model that further imposes orthogonal invariance of the summand distributions when viewed as symmetric tensors.
For the matrix case $p = 2$ earlier work identified the transition to Gaussian behavior at $r\asymp n^3$ for Gaussian summand vectors~\cite{BDER-2016-WignerWishartDetection}, and extended this result up to logarithmic factors to vectors with log-concave entry distributions~\cite{BG-2018-EntropicCLTWishart}.
These works show that, in our notation, $r_{\comp} \sim r_{\stat}$ on the polynomial scale when $p = 2$, so the above speculations of \cite{KMW-2024-TensorCumulantsInvariantInference} only apply to the case $p \geq 3$.

As we describe below, our main results show that, once $p \geq 3$, certain high-dimensional CLTs for tensors indeed \emph{do} have statistical-to-computational gaps, though \emph{not} the particular example proposed by \cite{KMW-2024-TensorCumulantsInvariantInference}.
We leave it as an intriguing direction for future research to understand what other limit theorems of classical probability theory exhibit such phenomena and thus admit substantially different computationally-bounded counterparts.
Tempting candidates in our opinion include Poisson approximations and limit theorems, mixing time and cutoff results, or asymptotic independence and the so-called propagation of chaos.

\subsection{Summary of contributions}
Our first goal in this paper is to clarify the above situation.
We show that the suggestion of \cite{KMW-2024-TensorCumulantsInvariantInference} is \emph{incorrect} and that the above high-dimensional CLT with Gaussian vectors $v$ does \emph{not} have a statistical-to-computational gap.
The reason for this is that the bound of \cite{Mikulincer-2020-CLTWishartTensors} on distributional distances is \emph{not} tight; we prove a sharpened bound in the TV distance (Theorem~\ref{thm:full-tensor-tv}) whose statistical threshold matches the computational one (Theorem~\ref{thm:full-tensor-polynomial-tests}).
A corresponding bound in Wasserstein distance (Lemma~\ref{lem:stein-error-distance} and Proposition~\ref{prop:tv-block-reduction}) also gives a direct sharpening of the results of \cite{Mikulincer-2020-CLTWishartTensors}.
As mentioned above, due to technical limitations \cite{KMW-2024-TensorCumulantsInvariantInference} did not quite prove a computational lower bound in the same model that \cite{Mikulincer-2020-CLTWishartTensors} studied, while here we also resolve this discrepancy and work with a single model (the more natural one of \cite{Mikulincer-2020-CLTWishartTensors}) for both of the results referenced above.

\begin{figure}
    \centering
    \begin{tikzpicture}[x=1.2cm, y=1cm, font=\small]
        \begin{scope}
            \node[anchor=west, font=\small\bfseries] at (0,1.2)
                {Undiluted tensor CLT: $c=1$};
            \fill[cgreen!25] (0,-0.22) rectangle (7.8,0.82);
            \fill[cred!25] (7.8,-0.22) rectangle (10.5,0.82);
            \node[text=cgreen!55!black] at (3.9,0.4) {easy};
            \node[text=cred!75!black] at (9.15,0.4) {impossible};
            \draw[thick] (0,0) -- (10.5,0);
            \draw[thick] (0,-0.14) -- (0,0.14);
            \fill (10.5,-0.12) -- (10.78,0) -- (10.5,0.12) -- cycle;
            \fill (7.8,0) circle (1.8pt);
            \node[below] at (7.8,-0.26)
                {$r_{\stat}, r_{\comp}\asymp n^3$};
        \end{scope}
        \begin{scope}[yshift=-2.4cm]
            \node[anchor=west, font=\small\bfseries] at (0,1.2)
                {Diluted tensor CLT: $c=n^{-2/3}$};
            \fill[cgreen!25] (0,-0.22) rectangle (2.8,0.82);
            \fill[yellow!35] (2.8,-0.22) rectangle (5.8,0.82);
            \fill[cred!25] (5.8,-0.22) rectangle (10.5,0.82);
            \node[text=cgreen!55!black] at (1.4,0.4) {easy};
            \node[text=yellow!45!black] at (4.3,0.4) {hard};
            \node[text=cred!75!black] at (8.15,0.4) {impossible};
            \draw[thick] (0,0) -- (10.5,0);
            \draw[thick] (0,-0.14) -- (0,0.14);
            \fill (10.5,-0.12) -- (10.78,0) -- (10.5,0.12) -- cycle;
            \fill (2.8,0) circle (1.8pt);
            \fill (5.8,0) circle (1.8pt);
            \node[below] at (2.8,-0.26) {$r_{\comp}\lesssim n^{1/3}$};
            \node[below] at (5.8,-0.26) {$r_{\stat}\gtrsim n^{2/3}$};
            \node at (5.25,-1.6) {$(r = \text{number of CLT summands})$};
        \end{scope}
    \end{tikzpicture}
    \caption{We illustrate regimes for the Wigner--Wishart CLT testing problem (Definition~\ref{def:wigner-wishart}) for 3-tensors with underlying entry distribution $\nu = \mathcal{N}(0, 1)$.
    Without dilution of the summands, the statistical and computational thresholds coincide at $r \asymp n^3$ (Theorems~\ref{thm:full-tensor-tv} and~\ref{thm:full-tensor-polynomial-tests}).
    With dilution rate $c=n^{-2/3}$, testing is hard for polynomials of degree $n^{\Omega(1)}$ once $r \gg n^{1/3 + \delta}$ for any fixed $\delta > 0$ (Theorem~\ref{thm:gaussian-dilution}), while brute-force search succeeds for $r\ll c^2n^2=n^{2/3}$ (Theorem~\ref{thm:scan-detection}), giving a gap (Corollary~\ref{cor:dilution-gap}).
    Theorems~\ref{thm:diluted-polynomial-tests} and~\ref{thm:diluted-tensor-tv} give matching efficient tests and TV distance bounds showing that these are indeed the computational and statistical thresholds, respectively.}
    \label{fig:clt-hardness}
\end{figure}
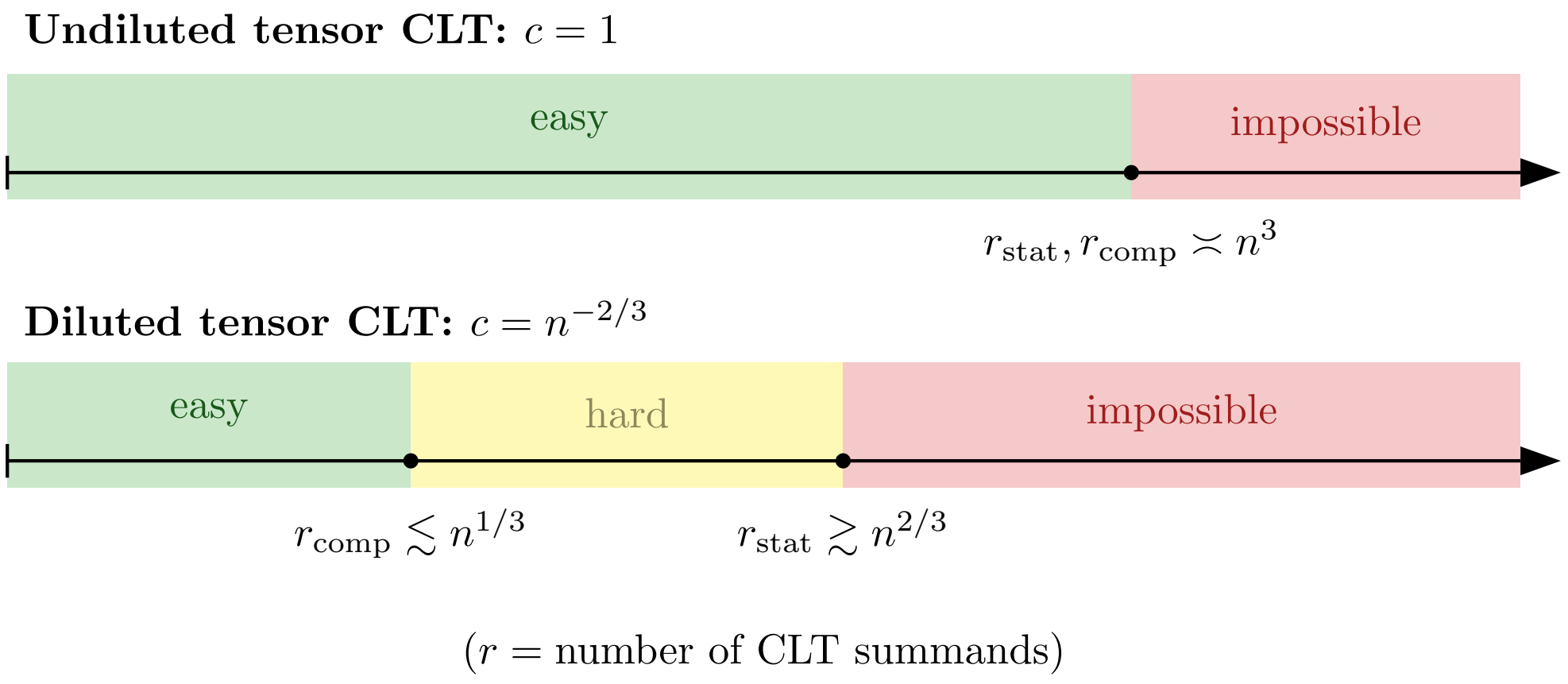

Our second goal is to show (to the best of our knowledge for the first time) that statistical-to-computational gaps \emph{do} occur for other high-dimensional CLTs, and to study what models give rise to them.
We show that the above tensor model has statistical-to-computational gaps (Corollary~\ref{cor:dilution-gap}) following a \emph{dilution} of the summand distribution $\mu$, replacing $\mu$ with $(c \cdot \mu) * (\sqrt{1 - c^2} \cdot \mathcal{N})$ for a small $c$ (equivalently, applying the Ornstein-Uhlenbeck process to $\mu$ for a long time $\log \frac{1}{c}$) and making the summands closer to Gaussian.
This operation makes testing harder and makes both $r_{\comp}$ and $r_{\stat}$ smaller; we show that, provided $c = c(n)$ is small enough, these two thresholds in fact separate on the polynomial scale in $n$.
The statistical test we use is a simple brute-force search over the typical set of the summand distribution, and applies to more general CLT models (Theorem~\ref{thm:scan-detection}).
We illustrate the case $p = 3$ of these results for a sufficiently small $c$ in Figure~\ref{fig:clt-hardness}.

For our analysis of computational CLTs via low-degree polynomials, we develop new tools for bounding \emph{low-degree $\chi^2$-divergences} (essentially the same as the so-called low-degree advantages; see Definition~\ref{def:low-deg-div}) through the cumulants of the distributions involved, inspired by previous works such as \cite{SW-2020-LowDegreeEstimation,RSWY-2022-CountCommunitiesLowDegree,LG-2024-GraphonEstimationLowDegree,SW-2025-SharpLowDegreeEstimation,EGV-2026-PermutationCumulants} identifying cumulant-like combinatorics in calculations with low-degree polynomial algorithms.
To more thoroughly demonstrate these tools, we further study computational CLTs in versions of matrix and tensor CLTs with \emph{masking}, where we only observe a subset of entries (Theorem~\ref{thm:masked-tensor-model}), akin to the masked Wishart model of \cite{BBH-2021-DeFinettiWishartMatrices} and related random geometric graph models \cite{BB-2024-FourierRandomGeometricGraphs,GPS-2026-BipartiteLatentSpaceGraphs}.

\subsection{AI methodology}
\label{sec:ai}

We used OpenAI's Codex, mostly with various tiers of the GPT-6 Astra model, to help develop and check proofs, revise prose, and edit the LaTeX source of this paper.
The biggest mathematical role played by these tools was in discovering the combination of techniques used in the rather involved proof of Theorem~\ref{thm:full-tensor-tv}, of which we present an edited version here.

\section{Main results}

\subsection{Formal setup of testing for CLTs}

We begin by formally defining the CLT model suggested above, viewed as a hypothesis testing problem.

\begin{definition}
    \label{def:clt-testing}
    Let $N = N(n), r = r(n) \in \mathbb{N}$ and let $\mu_n$ be a sequence of probability measures on $\mathbb{R}^N$ called the \emph{summand distributions}.
    Then, the \emph{CLT testing problem} with these parameters consists of the pair of probability measures:
    \begin{enumerate}
        \item $\mathbb{Q}_n = \mathcal{N} = \mathcal{N}(0, I_N)$.
        \item $\mathbb{P}_n = r^{-1/2} \cdot \mu_n^{*r}$. Equivalently, to sample $Y \sim \mathbb{P}_n$, we draw $X_1, \dots, X_r \sim \mu_n$ i.i.d.\ and observe $Y = r^{-1/2}\sum_{i=1}^r X_i$.
    \end{enumerate}
    We always use the letter $Y$ for observations from the measures $\mathbb{Q}_n$ and $\mathbb{P}_n$, and $X_i$ for the summands of $Y \sim \mathbb{P}_n$.
    We refer to $\Q_n$ as the \emph{null distribution} or \emph{model} and $\P_n$ as the \emph{alternative distribution} or \emph{model}.
\end{definition}
\noindent
We next introduce the notions of asymptotic success that we use for the hypothesis-testing problem above.
\begin{definition}[Strong detection]
    We say that a sequence of tests
\(
    t_n:\mathbb R^N\to\{0,1\}
\)
achieves \emph{strong detection} between $\mathbb Q_n$ and $\mathbb P_n$  if
\[
    \lim_{n \to \infty}\mathbb Q_n \left[t_n(Y)=0\right]
    =
    \lim_{n \to \infty}\mathbb P_n \left[t_n(Y)=1\right]
    =
    1.
\]
\end{definition}
\noindent
The following more quantitative version is the one we will work with in our main results on low-degree polynomial algorithms.
\begin{definition}[Strong separation]
    We say that a sequence of statistics $f_n:\mathbb R^N\to\mathbb R$ achieves \emph{strong separation} between $\mathbb Q_n$ and $\mathbb P_n$ if, as $n \to \infty$,
    \[ \Ex_{Y \sim \P_n} f_n(Y) - \Ex_{Y \sim \Q_n} f_n(Y) = \omega\left(\sqrt{\Varx_{Y \sim \Q_n} f_n(Y)} + \sqrt{\Varx_{Y \sim \P_n} f_n(Y)} \right). \]
\end{definition}

If $f_n(Y)$ achieves strong separation, then it is easy to check that setting $t_n(Y) \colonequals \bm 1\{f_n(Y) > \tau\}$ for a suitable threshold $\tau$ achieves strong detection, in essentially the same runtime as it takes to compute $f_n$.
Conversely, if $t_n(Y)$ achieve strong detection then these functions themselves also achieve strong separation.
Thus, the two notions are equivalent over arbitrary $f_n$ and $t_n$; however, as discussed in, e.g., the general survey \cite{Wein-2025-LowDegreeSurvey}, below we will consider strong separation by low-degree polynomial functions, which does not in general have a direct correspondence with strong detection by a special class of $t_n$.
See that survey for further discussion and justification of these notions, as well as reasoning behind taking strong separation by low-degree polynomials as a reasonable proxy for arbitrary computationally-efficient hypothesis tests.
This is the framework in which we formulate our computational lower bounds below.

\subsection{Wigner--Wishart tensor testing: No statistical-to-computational gap}
\label{sec:tensor-clt}

We now specialize the preceding framework to tensor-valued random variables.

\begin{definition}[Wigner--Wishart tensor testing problem]
\label{def:wigner-wishart}
Fix an integer $p\geq 2$ and a probability measure $\nu$ on $\mathbb R$
satisfying $\E_{v\sim\nu}[v]=0$ and $\E_{v\sim\nu}[v^2]=1$.
For each $n$, let
\[
    N \colonequals \binom{n}{p},
\]
and identify the coordinates of $\mathbb R^N$ with the $p$-element
subsets of $[n]\colonequals \{1,\ldots,n\}$, which we denote by $\binom{[n]}{p}$. Let $v \sim \nu^{\otimes n}$, i.e., $v_1, \dots, v_n \sim \nu$ drawn i.i.d., and denote the off-diagonal entries of the $p$-tensor $v^{\otimes p}$ by
\[
    X(v)
    \colonequals
    \left(
        \prod_{i\in S}v_i
    \right)_{S\in\binom{[n]}p}
    \in\mathbb R^N.
\]
Denote the law of this random vector by
\(
    \mu_n\colonequals \Law(X(v)).
\) For a sequence $r=r(n)\in\mathbb N$, the
\emph{Wigner--Wishart tensor testing problem} associated with the \emph{parameters} $(p,\nu,r)$
is the CLT testing problem of Definition~\ref{def:clt-testing} with
summand distribution $\mu_n$. Thus, the distributions involved are
\[
    \mathbb Q_n=\mathcal N,\qquad
    \mathbb P_n
    =
    r^{-1/2}\cdot\mu_n^{*r}.
\]
Equivalently, drawing $Y \sim \mathbb P_n$ we observe
\[
    Y
    =
    r^{-1/2}
    \sum_{a=1}^r X\bigl(v^{(a)}\bigr),
\]
where $v^{(1)},\ldots,v^{(r)}$ are independent copies of $v$.
\end{definition}
\noindent
We note that $\mu_n$ defined in this way is centered and isotropic.
Throughout we always use the letter $v$ for the vectors involved in the summand distributions of such constructions.

For $\nu$ with a finite third moment, we write $\tau\colonequals \Ex_{v\sim\nu}[v^3]$ and define the threshold scale
\[
    r_*(n) = r_*(n; p, \nu) \colonequals
    \begin{cases}
        n^{3p/2}, & p \text{ even},\\
        n^{(3p-1)/2}, & p \text{ odd and } \tau\neq0,\\
        n^p, & p \text{ odd and } \tau=0.
    \end{cases}
\]
As we will see, for Gaussian vectors $v^{(a)}$ the computational and statistical thresholds coincide at this scale, so we may view $r_{\comp}=r_{\stat}=r_*$.
Our first theorem gives convergence in TV distance above this scale for Gaussian entry distribution $\nu=\mathcal N(0,1)$.
The proof also gives a bound in quadratic Wasserstein distance with the same rates for $r\ge C_p n^p$ (Lemma~\ref{lem:stein-error-distance} and Proposition~\ref{prop:tv-block-reduction}) for a sufficiently large constant $C_p > 0$.
This may be compared directly with \cite[Theorem~3]{Mikulincer-2020-CLTWishartTensors}, improving its scaling condition $r\gg n^{2p-1}$ for Wasserstein convergence when $p\ge3$.
\begin{theorem}[Total variation distance bound]
\label{thm:full-tensor-tv}
Consider the Wigner--Wishart tensor testing problem of Definition~\ref{def:wigner-wishart}, with $\nu=\mathcal N(0,1)$ and $p\ge2$ fixed as $n\to\infty$.
There is a constant $C=C(p)$ such that
\[
    d_{\mathrm{TV}}(\mathbb P_n,\mathbb Q_n)
    \le C \cdot
    \begin{cases}
        \displaystyle \sqrt{\frac{n^{3p/2}}{r}},
        & p \text{ even},\\[1em]
        \displaystyle \frac{n^p}{r},
        & p \text{ odd}.
    \end{cases}
\]
In particular, for any $p \geq 2$, if $r=\omega(r_*(n))$, then $d_{\mathrm{TV}}(\mathbb P_n,\mathbb Q_n)\to0$.
\end{theorem}

\begin{remark}[Non-Gaussian vectors]
We expect that our argument for Theorem~\ref{thm:full-tensor-tv} can be adapted, with additional technical work, to at least two further classes of input laws $\nu$ with mean zero and variance one.
The first consists of laws with density proportional to $\exp(-V(v))$, where $V\in C^\infty(\mathbb R)$, $0<c_0\le V''(v)\le C_0<\infty$, and $V^{(j)}$ is bounded for every $j\ge3$ (a strong combination of smoothness and strict log-concavity).
The second consists of mean zero and variance one rescalings of $\operatorname{Beta}(a,b)$ laws with fixed $a,b\ge1$.
For both classes, the resulting bound should give convergence in TV distance when $r=\omega(r_*(n))$ for the same threshold scale.
\end{remark}
\noindent
The next theorem gives polynomial-time tests below $r_*(n)$, showing that the threshold in Theorem~\ref{thm:full-tensor-tv} is tight.
\begin{theorem}[Polynomial tests below the threshold]
\label{thm:full-tensor-polynomial-tests}
Consider the Wigner--Wishart tensor testing problem of Definition~\ref{def:wigner-wishart}, with $p\ge2$ and $\nu$ fixed as $n\to\infty$.
Assume that $\Ex_{v\sim\nu}[|v|^6]<\infty$.
If $r=o(r_*(n))$, there is an explicit sequence of polynomials $f_n$ of degree at most four achieving strong separation between $\mathbb P_n$ and $\mathbb Q_n$.
Each $f_n$ can be evaluated in $O(n^{2p})$ arithmetic operations.
Thresholding these statistics achieves strong detection, and $d_{\mathrm{TV}}(\mathbb P_n,\mathbb Q_n)\to1$.
\end{theorem}
\noindent
The proof is in Section~\ref{sec:pf:thm:full-tensor-polynomial-tests}.
Thus, for Gaussian vectors $v^{(a)}$, the statistical and computational thresholds coincide at $r_*(n)$.

For the sake of completeness, we also give a low-degree lower bound under only the weaker moment condition of Definition~\ref{def:bounded-moments} on $\nu$.

\begin{theorem}
\label{thm:low-degree-1}
Let $(\mathbb Q_n,\mathbb P_n)$ be the Wigner--Wishart tensor testing
problem of Definition~\ref{def:wigner-wishart}, with parameters
$p$, $\nu$, and $r=r(n)$. Assume that $\nu$ has $(K,C)$-bounded moments (Definition~\ref{def:bounded-moments}).
There is a constant $\gamma_{p,\nu}>0$ such that the following holds for a sufficiently small constant $\kappa=\kappa(p,K,C)>0$.
If we have
\[
    r=\omega \left(r_*(n)\right), \quad D \le\kappa\left(
        \frac{r}{r_*(n)}
    \right)^{1/\gamma_{p,\nu}},
\]
then no sequence of polynomials of degree at most $D$ strongly separates $\mathbb P_n$ and $\mathbb Q_n$.
\end{theorem}
\noindent
Corresponding upper and lower bounds on low-degree $\chi^2$-divergences are stated in Theorem~\ref{thm:low-degree-1-chi-squared}, and an explicit choice of $\gamma_{p,\nu}$ is given in Definition~\ref{defin:q,s}.

\subsection{Diluted models and statistical-to-computational gaps}
\label{sec:dilution-results}

We next turn to the procedure that we show can generate statistical-to-computational gaps.
\begin{definition}[Gaussian dilution]
Let $\mu$ be a probability measure on $\mathbb R^N$ and let
$c=c(n)\in[0,1]$. The \emph{Gaussian dilution} of $\mu$
at level $c$, denoted by $G_c(\mu)$, is defined as
\[
    G_c(\mu)
    \colonequals
    \Law \left(
        cY+\sqrt{1-c^2}\,G
    \right) = (c \cdot \mu) * (\sqrt{1 - c^2} \cdot \mathcal{N}),
\]
where in the first expression
$Y\sim\mu$ and $G\sim\mathcal N(0,I_N)$ are independent.
\end{definition}
\noindent
The parameter $c$ therefore interpolates between the standard Gaussian measure and
the original distribution: $G_0(\mu)=\mathcal N$ and
$G_1(\mu)=\mu$.

\begin{remark}
    The standard Gaussian measure is invariant under Gaussian dilution, and Gaussian dilution commutes with normalized convolution.
    In the setting of the CLT testing model, Gaussian dilution acts as
\begin{align*}
  G_c(\mathbb Q_n) &= \mathbb Q_n, \\
  G_c(\mathbb P_n)
    &=
    G_c \left(
        r^{-1/2}\cdot\mu_n^{*r}
    \right)
    =
    r^{-1/2}\cdot G_c(\mu_n)^{*r}.
\end{align*}
\end{remark}

We now consider applying Gaussian dilution to the Wigner--Wishart tensor model and its effect on the computational and statistical thresholds.
\begin{theorem}[Low-degree lower bound under dilution]
\label{thm:gaussian-dilution}
Let $(\mathbb Q_n,\mathbb P_n)$ be the Wigner--Wishart tensor testing
problem of Definition~\ref{def:wigner-wishart}.
Assume that $\nu$ has $(K,C)$-bounded moments (Definition~\ref{def:bounded-moments}).
For $c=c(n)\in(0,1]$, define
\[
    r_{\comp}(n; c)
    \colonequals
    \max\left\{
        c^6r_*(n),
        c^4n^p
    \right\}.
\]
Let $\gamma_{p,\nu}$ be as in Theorem~\ref{thm:low-degree-1}.
For a sufficiently small constant $\kappa=\kappa(p,K,C)>0$, if we have
\[
    r=\omega\bigl(r_{\comp}(n; c)\bigr), \quad
    D\le\kappa
    \left(\frac{r}{r_{\comp}(n; c)}\right)^{1/\gamma_{p,\nu}},
\]
then no sequence of polynomials of degree at most $D$ achieves strong separation between $G_c(\mathbb P_n)$ and $G_c(\mathbb Q_n) = \mathbb Q_n$, or equivalently in the CLT testing problem with summand distribution $G_c(\mu_n)$.
\end{theorem}
\noindent
The corresponding analyses of low-degree $\chi^2$-divergences are given in
Theorem~\ref{thm:gaussian-dilution-chi-squared}.
The next theorem gives matching polynomial tests below this scale.
\begin{theorem}[Polynomial tests under dilution]
\label{thm:diluted-polynomial-tests}
Consider the Wigner--Wishart tensor testing problem of Definition~\ref{def:wigner-wishart}, with $p\ge2$ and $\nu$ fixed as $n\to\infty$.
Assume that $\nu$ has $(K,C)$-bounded moments (Definition~\ref{def:bounded-moments}).
Let $c=c(n)\in(0,1]$ and let $r_{\comp}(n;c)$ be as in Theorem~\ref{thm:gaussian-dilution}.
If $r=o(r_{\comp}(n;c))$, there is an explicit sequence of polynomials $f_n$ of degree at most four achieving strong separation between $G_c(\mathbb P_n)$ and $\mathbb Q_n$.
Each $f_n$ can be evaluated in $O(n^{2p})$ arithmetic operations.
Thresholding these statistics achieves strong detection, and $d_{\mathrm{TV}}(G_c(\mathbb P_n),\mathbb Q_n)\to1$.
\end{theorem}
\noindent
The proof is in Section~\ref{sec:pf:thm:diluted-polynomial-tests}.

For sufficiently small $c$, a brute-force search can succeed beyond the range of these polynomial tests by finding a single summand $X_i$ unusually well aligned with the observation $Y \sim G_c(\mathbb{P}_n)$.
The following theorem gives a general-purpose formulation of this strategy and applies
to any centered and isotropic summand distribution.

\begin{theorem}[Brute-force summand search]
\label{thm:scan-detection}
Consider the CLT testing problem of Definition~\ref{def:clt-testing}, with
centered and isotropic summand distributions $\mu_n$ on $\mathbb R^N$,
and let $c=c(n)\in(0,1]$.
Let $\mathcal T_n$ be a deterministic, nonempty finite set of unit vectors
in $\mathbb R^N$, and write $M_n\colonequals|\mathcal T_n|$.
Suppose that, for some $a_n>0$,
\[
    \Prx_{X\sim\mu_n}\left[
        \max_{v\in\mathcal T_n}\langle v,X\rangle\ge a_n
    \right]\longrightarrow1.
\]
If we have
\[
    r(1+\log M_n)=o(c^2a_n^2),
\]
then the test $t_n$ defined by
\[
    T_n(y)\colonequals\max_{v\in\mathcal T_n}\langle v,y\rangle,
    \qquad
    t_n(y)\colonequals\bm 1\left\{T_n(y)>\frac{ca_n}{2\sqrt r}\right\}
\]
achieves strong detection between $G_c(\mathbb P_n)$ and $\mathbb Q_n$.

In particular, suppose that $\mathcal S_n\subseteq\mathbb R^N\setminus\{0\}$
is a measurable set with $\mu_n(\mathcal S_n)\to1$ and
$\|x\|_2\ge b_n>0$ for every $x\in\mathcal S_n$.
Then, either of the following choices supplies a valid set of test directions $\mathcal T_n$:
\begin{enumerate}
    \item \emph{Discrete summands:} If $\mathcal S_n$ is finite, take
    $\mathcal T_n=\{x/\|x\|_2:x\in\mathcal S_n\}$ and $a_n=b_n$.
    It suffices that $r(1+\log|\mathcal S_n|)=o(c^2b_n^2)$.
    \item \emph{Net of typical set:} Take $\mathcal T_n$ to be a finite $\varepsilon_n$-net of
    $\{x/\|x\|_2:x\in\mathcal S_n\}$ on the unit sphere, where
    $0<\varepsilon_n<\sqrt2$, and take
    $a_n=(1-\varepsilon_n^2/2)b_n$.
\end{enumerate}
\end{theorem}

\begin{corollary}[Statistical strong detection under dilution]
\label{cor:diluted-tensor-search}
Consider the Wigner--Wishart tensor testing problem of Definition~\ref{def:wigner-wishart}, with $p\ge2$ and $\nu$ fixed as $n\to\infty$.
Assume that $\nu$ has $(K,C)$-bounded moments (Definition~\ref{def:bounded-moments}).
Let $c=c(n)\in(0,1]$, let $r_{\comp}(n;c)$ be as in Theorem~\ref{thm:gaussian-dilution}, and define
\[
    r_{\stat}(n;c)=r_{\stat}(n;c,p,\nu)
    \colonequals\max\{c^2n^{p-1},r_{\comp}(n;c)\}.
\]
If $r=o(r_{\stat}(n;c))$, there is a sequence of tests achieving strong detection between $G_c(\mathbb P_n)$ and $\mathbb Q_n$, using a brute-force search over $\exp(O_p(n))$ directions when $c^2n^{p-1}\ge r_{\comp}(n;c)$, and the polynomial tests of Theorem~\ref{thm:diluted-polynomial-tests} otherwise.
In particular, $d_{\mathrm{TV}}(G_c(\mathbb P_n),\mathbb Q_n)\to1$.
\end{corollary}
\noindent
The proof is in Section~\ref{sec:pf:thm:scan-detection}.

For tensor summands with Gaussian underlying vectors, the following convergence result complements these tests, showing that this $r_{\stat}(n; c)$ is indeed the statistical threshold on the polynomial scale.
\begin{theorem}[Total variation distance bound under dilution]
\label{thm:diluted-tensor-tv}
Consider the Wigner--Wishart tensor testing problem of Definition~\ref{def:wigner-wishart}, with $\nu=\mathcal N(0,1)$ and $p\ge2$ fixed as $n\to\infty$.
Let $c=c(n)\in(0,1]$ and let $r_{\stat}(n;c)$ be as in Corollary~\ref{cor:diluted-tensor-search}.
For any fixed $\delta>0$, if
\[
    r\ge n^\delta r_{\stat}(n;c),
\]
then $d_{\mathrm{TV}}(G_c(\mathbb P_n),\mathbb Q_n)\to0$.
\end{theorem}
\noindent
The proof in Section~\ref{sec:pf:thm:diluted-tensor-tv} is a simplified variant of the proof of the undiluted case in Theorem~\ref{thm:full-tensor-tv}.

We now compare the statistical and computational threshold scales.
For general $\nu$ satisfying the moment assumption above, $r_{\stat}(n;c)$ is a lower bound on the statistical threshold.
Theorems~\ref{thm:gaussian-dilution} and~\ref{thm:diluted-polynomial-tests} identify $r_{\comp}(n;c)$ as the computational threshold in the low-degree sense used here.
For Gaussian $\nu$, Theorem~\ref{thm:diluted-tensor-tv} and Corollary~\ref{cor:diluted-tensor-search} identify the statistical threshold on the polynomial scale as $r_{\stat}(n;c)$.
These results give a statistical-to-computational gap whenever
$r_{\comp}(n; c)=o(r_{\stat}(n; c))$ and $r_{\stat}(n; c)\to\infty$.
Indeed,
\[
    \frac{r_{\comp}(n; c)}{r_{\stat}(n; c)}
    =\min\left\{1,\max\left\{
        \frac{c^4r_*(n)}{n^{p-1}},
        c^2n
    \right\}\right\}.
\]
Define
\[
    \alpha_* = \alpha_*(p,\nu) \colonequals
    \begin{cases}
        (p+2)/8, & p\text{ even},\\
        (p+1)/8, & p\text{ odd and }\tau\ne0,\\
        1/2, & p\text{ odd and }\tau=0.
    \end{cases}
\]
Then the required range of the dilution parameter $c$ to observe a gap is
\[
    n^{-(p-1)/2}\ll c\ll n^{-\alpha_*}.
\]
This range is nonempty for every $p\ge3$; for $p=2$ (the matrix case), the two bounds on $c$ are incompatible, so this comparison does not give a gap.
We summarize this state of affairs below.

\begin{corollary}[Statistical-computational gaps under dilution]
\label{cor:dilution-gap}
Consider the tensor model of Definition~\ref{def:wigner-wishart}, with
fixed $p\ge3$ and a fixed law $\nu$ having $(K,C)$-bounded moments (Definition~\ref{def:bounded-moments}).
Let $\gamma_{p,\nu}$ be as in Theorem~\ref{thm:low-degree-1}.
Choose a constant $\alpha$ satisfying
\[
    \alpha_*<\alpha<\frac{p-1}{2},
\]
set $c=n^{-\alpha}$, and write
\[
    r_{\comp}(n;c)=n^\lambda.
\]
The interval
\[
    \max\{0,\lambda\}<\beta<p-1-2\alpha
\]
is nonempty. For any fixed $\beta$ in this interval and
$r=r(n) \colonequals \lfloor n^\beta\rfloor$, we have:
\begin{enumerate}
    \item A brute-force search over $\exp(O_p(n))$ directions achieves strong detection
    between $G_c(\mathbb P_n)$ and $\mathbb Q_n$.
    In particular,
    $d_{\mathrm{TV}}(G_c(\mathbb P_n),\mathbb Q_n)\to1$.
    \item For a sufficiently small constant $\kappa>0$ (not depending on $n$), no sequence of
    polynomials of degree at most
    \[
        D=\kappa \cdot n^{(\beta-\lambda)/\gamma_{p,\nu}}
    \]
    achieves strong separation between $G_c(\mathbb P_n)$ and $\mathbb Q_n$.
\end{enumerate}
\end{corollary}
For a concrete example, take cubic tensors ($p=3$) with Gaussian coordinates and $c=n^{-2/3}$, as in Figure~\ref{fig:clt-hardness}.
Then $r_{\comp}(n; c)=n^{1/3}$ and $r_{\stat}(n; c)=n^{2/3}$, so we may choose $r=\lfloor n^{1/2}\rfloor$ in the gap.
Brute-force search over a net of $\exp(O(n))$ directions achieves strong detection, while no polynomials of degree at most $n^\delta$ achieve strong separation for some constant $\delta>0$.

\begin{remark}[Other summand distributions]
We expect the same strategy to give gaps beyond tensor models as well.
For example, consider a centered isotropic law whose typical summands
have norm of order $\sqrt N$ and whose normalized directions have a
Euclidean $1$-net of size
$\exp(H_n)$, with $H_n=o(N)$.
The brute-force search then succeeds when $r\ll c^2N/(1+H_n)$.
If the first cumulants of $\mu_n$ not matching those of $\mathcal{N}$ occur at a fixed degree $q\ge3$, write
$A_n\colonequals\sum_{|\alpha|=q}\kappa_\alpha(\mu_n)^2/\alpha!$.
The weighted sum of their squares after convolution and dilution is
$c^{2q}r^{2-q}A_n$, suggesting the scale
$(c^{2q}A_n)^{1/(q-2)}$ for the computational threshold, which decreases faster with $c$ than the
above scale on which the statistical brute-force search succeeds.
However, uniform control of the higher cumulants and their decompositions, as in
Lemma~\ref{lem:clt-cumulant-support-bound} for our tensor models, would be needed to turn this
comparison into a rigorous gap for growing degrees.
\end{remark}

\subsection{Masked models}
\label{sec:masked-models}

We now turn to an equally natural procedure for making a CLT testing problem harder, but one which turns out to behave quite differently from dilution.

\begin{definition}[Masking]
\label{defin:masking}
Let $\mathbb P$ be a probability measure on $\mathbb R^N$, and let
$\Omega\subseteq[N]$ be a set of observed coordinates. Define the coordinate
projection
\[
    \Pi_\Omega:\mathbb R^N\to\mathbb R^{|\Omega|},
    \qquad
    \Pi_\Omega(x)\colonequals (x_e)_{e\in\Omega}.
\]
The \emph{masked distribution} of $\mathbb P$ with respect to $\Omega$,
denoted by $M_\Omega(\mathbb P)$, is the law of $\Pi_\Omega(Y)$
when $Y\sim\mathbb P$.
\end{definition}

\begin{remark}
In the CLT testing problem of Definition~\ref{def:clt-testing}, masking
commutes with normalized convolution. In particular,
\begin{align*}
      M_\Omega(\mathbb Q_n)
  &=
    \mathcal N(0,I_{|\Omega|}), \\
    M_\Omega(\mathbb P_n)
    &=
    M_\Omega \left(
        r^{-1/2}\cdot\mu_n^{*r}
    \right)
    =
      r^{-1/2}\cdot M_\Omega(\mu_n)^{*r}.
      \end{align*}
\end{remark}

We next introduce the hypergraph notation needed to describe thresholds we derive for masked tensor CLT models.
For a multi-index $\alpha=(\alpha_1,\ldots,\alpha_N)\in\mathbb Z_{\geq 0}^N$ and $x=(x_1,\ldots,x_N)\in\mathbb R^N$, we write $|\alpha|\colonequals\sum_{i=1}^N\alpha_i$, $\alpha!\colonequals\prod_{i=1}^N\alpha_i!$, and $x^\alpha\colonequals\prod_{i=1}^N x_i^{\alpha_i}$.
A multi-index $\alpha\in\mathbb Z_{\geq 0}^{N}$ may, under the identification of $[N]$ with $\binom{[n]}{p}$,
be identified with a $p$-uniform multihypergraph: for each
$e\in\binom{[n]}p$, the hyperedge $e$ appears with multiplicity
$\alpha_e$. We denote this multihypergraph by $H_\alpha$. Its number of
edges, counted with multiplicity, is $|E(H_\alpha)|=|\alpha| = \sum_{e}\alpha_e$.
We use the term \emph{shape} for the isomorphism class of such $p$-uniform multihypergraphs.

Since the coordinates $v_1,\ldots,v_n$ are independent and centered, basic properties of cumulants given in
Propositions~\ref{prop:moment-cumulant-identities} and~\ref{prop:cumulant-independence} imply that a cumulant $\kappa_\alpha(\mu_n)$ vanishes whenever $H_\alpha$ is disconnected or contains a vertex of degree one (see Section~\ref{sec:tv-cumulants}).
Let $q=4$ when $p$ is odd and $\tau=0$, and $q=3$ otherwise.
This is the first degree at which the cumulants of $\mu_n$ differ from those of $\mathcal N$.
Accordingly, we will consider
\[
    \mathcal A \coloneqq
    \left\{
        \alpha\in\mathbb Z_{\geq 0}^{\binom{[n]}{p}}:
        |\alpha|\geq q,\;
        H_\alpha \text{ is connected},\;
        \deg_{H_\alpha}(v)\geq 2
        \text{ for every }v\in V(H_\alpha)
    \right\}.
\]
We refer to a shape corresponding to an element of $\mathcal A$ as an \emph{admissible shape}.

For a shape $H$, let $f(H)$ denote the largest number of admissible
subshapes into which its edge multiset can be partitioned. More precisely,
\[
    f(H)
    \coloneqq
    \max\left\{
        k:
        E(H)
        =
        E(H_1)\sqcup\cdots\sqcup E(H_k),
        H_1,\ldots,H_k \text{ admissible}
    \right\},
\]
where the decomposition respects edge multiplicities. If no such
decomposition exists, we set $f(H)=0$.

\begin{definition}[Irreducible shape]
\label{def:irreducible}
An admissible shape $H$ is called \emph{irreducible} if $f(H)=1$.
Equivalently, $H$ is irreducible if its edge multiset cannot be partitioned into two or more nonempty admissible shapes.
\end{definition}

We now consider the effect of masking. Let $\Omega\subseteq\binom{[n]}p$
be a deterministic collection of observed coordinates. As in Definition~\ref{defin:masking}, we write $M_\Omega(\mathbb P_n)$ and
$M_\Omega(\mathbb Q_n)$ for the corresponding masked alternative and null distributions in a CLT testing problem.
For a $p$-uniform multihypergraph shape $H$, define
\[
    \rho(H,\Omega) \coloneqq
    \#\left\{\phi:V(H)\hookrightarrow[n]:\phi(e)\in\Omega \text{ for every }e\in E(H)
    \right\}.
\]
Thus, $\rho(H,\Omega)$ is the number of injective embeddings of $H$ into the hypergraph determined by the mask $\Omega$.
For a degree cutoff $D$, also define
\[
    \zeta_{\leq D}(r; \Omega) \coloneqq
    \max_{\substack{
        H\ {\rm irreducible}\\
        q\leq |E(H)|\leq D}}
    \rho(H,\Omega)^{1/|E(H)|} \cdot
    r^{2/|E(H)|-1}.
\]
We will see that $\zeta_{\leq D}(r; \Omega)$ measures the largest contribution per edge of any irreducible admissible configuration supported by the mask.

Our main result then gives a computational lower bound in terms of these quantities.
\begin{theorem}[Masked Wigner--Wishart tensors]
\label{thm:masked-tensor-model}
Let $(\mathbb Q_n,\mathbb P_n)$ be the Wigner--Wishart tensor testing
problem of Definition~\ref{def:wigner-wishart}.
Assume that $\nu$ has $(K,C)$-bounded moments (Definition~\ref{def:bounded-moments}).
For a sufficiently small constant $\kappa=\kappa(p,K,C)>0$, if we have
\[
    \zeta_{\leq D}(r; \Omega)=o(1), \qquad D
    \le\kappa
    \left(
        \frac{1}{\zeta_{\leq D}(r; \Omega)}
    \right)^{1/(2pK+p)},
\]
then no sequence of polynomials of degree at most $D$ achieves strong
separation between $M_\Omega(\mathbb P_n)$ and $M_\Omega(\mathbb Q_n)$.
\end{theorem}
\noindent
The corresponding bounds on low-degree $\chi^2$-divergences are stated in
Theorem~\ref{thm:masked-tensor-model-chi-squared}.

We remark that, while masking and especially masking random entries seems to be an operation analogous to dilution, regardless of how small a fraction of entries are revealed, perhaps surprisingly the previous proof technique from Section~\ref{sec:dilution-results} does \emph{not} apply to random masking to demonstrate a statistical-to-computational gap.
\begin{remark}[Bernoulli masking and brute-force search]
    Suppose, by analogy with dilution, that each entry is retained independently with probability $c=c(n)\in(0,1]$.
    As with dilution, decreasing $c$ makes the associated CLT testing problem harder.
    Write $s_{\mathrm{mask}}(n; c, p)\colonequals c n^{p-1}$.
    For fixed $p$ and $\nu$ under the moment assumptions above, when $s_{\mathrm{mask}}(n; c, p)\to\infty$, a typical masked summand has squared norm of order $c n^p$, with high probability over the mask.
    Theorem~\ref{thm:scan-detection}, applied without Gaussian dilution, therefore gives strong detection by brute-force search for $r=o(s_{\mathrm{mask}}(n; c, p))$ when $\nu$ has finite support.
    A net of typical summand directions with $\log|\mathcal T_n|=O(n\log(2/c))$ gives the condition $r=o(s_{\mathrm{mask}}(n; c, p)/\log(2/c))$ for general $\nu$.
    These brute-force search bounds are not vacuous only when $s_{\mathrm{mask}}(n; c, p)\to\infty$.

    To compare this with Theorem~\ref{thm:masked-tensor-model}, consider the irreducible shape $H$ consisting of two distinct $p$-edges intersecting in one vertex, each repeated twice.
    It has four edges counted with multiplicity, but requires only two distinct entries to be observed, so with high probability
    \[
        \rho(H,\Omega)\asymp c^2n^{2p-1},
        \qquad
        \zeta_{\leq D}(r; \Omega)\ge \rho(H,\Omega)^{1/4}r^{-1/2}
        \text{ for }D\ge4.
    \]
    Thus the hypothesis $\zeta_{\leq D}(r; \Omega)=o(1)$ already requires $r=\omega(c n^{p-1/2})=\omega( n^{1/2}s_{\mathrm{mask}}(n; c, p))$, beyond the range of the brute-force search bound.
    These two bounds therefore do \emph{not} yield a statistical-to-computational gap from Bernoulli masking alone.
\end{remark}

\section{Notation}

For a multi-index $\alpha=(\alpha_1,\ldots,\alpha_N)\in\mathbb Z_{\geq 0}^N$ and $x=(x_1,\ldots,x_N)\in\mathbb R^N$, we write $|\alpha|\colonequals\sum_{i=1}^N\alpha_i$, $\alpha!\colonequals\prod_{i=1}^N\alpha_i!$, and $x^\alpha\colonequals\prod_{i=1}^N x_i^{\alpha_i}$.
We denote by $\mathcal{N}$ the standard Gaussian measure on $\mathbb{R}^N$,
with the ambient dimension $N$ understood from context; that is, $\mathcal{N} \colonequals  \mathcal{N}(0,I_N)$.
For two probability measures $\mu$ and $\nu$ on $\mathbb R^N$,
we write $a\cdot\mu$ for the law of $aX$ when $X\sim\mu$,
and $\mu*\nu$ for their convolution, namely the law of the sum
of two independent random vectors with respective laws $\mu$
and $\nu$: $\mu * \nu \colonequals \Law(X+Y)$, where $X \sim \mu$ and $Y \sim \nu$ are independent.
More generally, for $r\in\mathbb{N}$, we write $\mu^{*r}$ for the $r$-fold convolution of $\mu$ with itself.
Equivalently, $\mu^{*r}$ is the law of $X_1+\cdots+X_r$ for independent $X_1,\ldots,X_r\sim\mu$.

\section{Low-degree lower bounds via generating functions}

In this section, we develop a toolkit based on moment and cumulant generating functions for bounding the low-degree $\chi^2$-divergences (Definition~\ref{def:low-deg-div}), which may be of independent interest.
We first review the general toolkit surrounding low-degree polynomial algorithms, and then give our new approach.

\subsection{Low-degree algorithms and divergences}
\label{sec:low-degree}

The basic principle of the analysis of low-degree polynomial algorithms for hypothesis testing is that bounded low-degree
$\chi^2$-divergence rules out strong separation by low-degree polynomials.
We refer to the surveys
\cite{kunisky2019notes,Wein-2025-LowDegreeSurvey}
for a more detailed treatment of the low-degree method.

\begin{definition}[Low-degree $\chi^2$-divergence]
\label{def:low-deg-div}
Let $\mathbb P$ and $\mathbb Q$ be probability measures on $\mathbb R^N$.
For $D\in\mathbb N$, the \emph{degree-$D$ $\chi^2$-divergence} of
$\mathbb P$ from $\mathbb Q$ is defined by
\begin{align*}
    \chi^2_{\leq D}(\P \mid \Q)
    &\colonequals \sup_{f \in \R[Y]_{\leq D} \setminus \{0\}} \frac{(\mathbb{E}_{\P}[f(Y)])^2}{\mathbb{E}_{\mathbb{Q}}[f(Y)^2]} - 1 \\
    &= \sup_{\substack{f \in \R[Y]_{\leq D} \setminus \{0\} \\ \mathbb{E}_{\Q} f(Y) = 0}}\frac{(\mathbb{E}_{\P}[f(Y)])^2}{\mathbb{E}_{\mathbb{Q}}[f(Y)^2]}.
\end{align*}
\end{definition}
\noindent
If the restriction to degree-$D$ polynomials is removed, the same
variational expression recovers the classical $\chi^2$-divergence
\[
    \chi^2(\mathbb P\mid\mathbb Q)
    \colonequals
    \E_{\mathbb Q}
    \left[
        \left(
            \frac{d\mathbb P}{d\mathbb Q}-1
        \right)^2
    \right],
\]
whenever $\mathbb P\ll\mathbb Q$ and
$d\mathbb P/d\mathbb Q\in L^2(\mathbb Q)$.
The quantity
\(
    \sqrt{1+\chi^2_{\leq D}(\mathbb P\mid\mathbb Q)}
\)
is often referred to as the \emph{low-degree advantage}. We use the
$\chi^2_{\leq D}$ notation to emphasize its close analogy with the
classical $\chi^2$-divergence.

The low-degree and classical $\chi^2$-divergences govern strong separation by low-degree polynomials and by arbitrary functions in the following way.
\begin{proposition}[Proposition 6.2 of \cite{BAHSWZ-2022-FranzParisiLowDegree}, Lemma 1.13 of \cite{kunisky2019notes}]
\label{prop:chi-squared-testing}
Consider a sequence of testing problems
$(\mathbb P_n,\mathbb Q_n)$ on $\mathbb R^{N(n)}$.
\begin{enumerate}
    \item If
    \(
        \chi^2_{\leq D(n)}
        (\mathbb P_n\mid\mathbb Q_n)
        =
        O(1),
    \)
    then no sequence of polynomials
    \(
        f_n\in\mathbb R[Y]_{\leq D(n)}
    \)
    achieves strong separation between $\mathbb P_n$ and $\mathbb Q_n$.

    \item
    If
    \(
        \chi^2(\mathbb P_n\mid\mathbb Q_n)=O(1),
    \)
    then no sequence of measurable functions
    \(
        f_n:\mathbb R^{N(n)}\to\mathbb R
    \)
    achieves strong separation, and no sequence of tests
    \(
        t_n:\mathbb R^{N(n)}\to\{0,1\}
    \)
    achieves strong detection.
\end{enumerate}
\end{proposition}
\noindent
We recall that neither implication admits a converse in general, and indeed these divergences can diverge while the underlying testing problem remains difficult; see, for example,
\cite{BAHSWZ-2022-FranzParisiLowDegree,
COGHWZ-2022-PhaseTransitionsGroupTesting,
DMW-2023-DenseSubhypergraphsLowDegree,
DDL-2023-GraphMatchingLowDegree}.

\subsection{Generating function representation}

We now develop a power series formalism for representing the low-degree
$\chi^2$-divergence in terms of moment generating functions (MGFs). The main
advantage of this representation is that a ratio of MGFs
can be written as the exponential of a difference of
cumulant generating functions (CGFs), allowing low-degree distinguishability to
be controlled directly through cumulants.

\begin{definition}
Let $\mu$ be a probability measure on $\mathbb R^N$ with
finite moments of all orders. Its formal MGF is
\[
    \phi_\mu(t)
    \coloneqq
    \sum_{\alpha\in\mathbb Z_{\geq 0}^N}
    \frac{m_\alpha(\mu)}{\alpha!}t^\alpha,
    \qquad
    m_\alpha(\mu)\colonequals \mathbb E_{X\sim\mu}[X^\alpha].
\]
Its formal CGF is
\[
    \psi_\mu(t)
    \colonequals
    \log\phi_\mu(t)
    =
    \sum_{\alpha\in\mathbb Z_{\geq 0}^N\setminus\{0\}}
    \frac{\kappa_\alpha(\mu)}{\alpha!}t^\alpha.
\]
The formal logarithm is well-defined because $\phi_\mu(0)=1$.
We interpret these series and their identities coefficientwise, without assuming convergence.
\end{definition}

\begin{proposition}[Multiplicativity of MGF]
Let $\mu$ and $\nu$ be probability measures on $\mathbb R^N$
with moments of all orders. Then,
\[
    \phi_{\mu*\nu}(t)=\phi_\mu(t)\phi_\nu(t),
    \qquad
    \psi_{\mu*\nu}(t)=\psi_\mu(t)+\psi_\nu(t).
\]
Consequently, for every integer $r\ge1$,
\[
    \phi_{\mu^{*r}}(t)=\phi_\mu(t)^r,
    \qquad
    \psi_{\mu^{*r}}(t)=r\psi_\mu(t).
\]
In particular, for every nonzero $\alpha$,
\[
    \kappa_\alpha(\mu*\nu)
    =\kappa_\alpha(\mu)+\kappa_\alpha(\nu).
\]
\end{proposition}

We recall that the multivariate central limit theorems we discussed above involve distributions of the form $r^{-1/2} \cdot \mu^{*r}$.
Below we derive how the cumulants of these behave.
\begin{proposition}[Cumulants under scaling]
\label{prop:scaling-cumulants}
For $c\in\mathbb R$ and a probability measure $\mu$ on $\mathbb R^N$,
\[
    \phi_{c\cdot\mu}(t)=\phi_\mu(ct),
    \qquad
    \psi_{c\cdot\mu}(t)=\psi_\mu(ct).
\]
Consequently, for every multi-index $\alpha$,
\[
    \kappa_\alpha(c\cdot\mu)
    =
    c^{|\alpha|}\kappa_\alpha(\mu).
\]
In particular,
\[
    \kappa_\alpha
    \left(
        r^{-1/2}\cdot\mu^{*r}
    \right)
    =
    r^{1-|\alpha|/2}\kappa_\alpha(\mu).
\]
\end{proposition}

\begin{definition}[Power series norms]
Let
\[
    F(t)
    =
    \sum_{\alpha\in\mathbb Z_{\geq 0}^N}
    \frac{c_\alpha}{\alpha!}t^\alpha
\]
be a formal power series with coefficients $c_\alpha\in\mathbb R$.
We define
\[
    \|F\|_{\mathcal P^2}^2
    \colonequals
    \sum_{\alpha\in\mathbb Z_{\geq 0}^N}
    \frac{|c_\alpha|^2}{\alpha!},
\]
For $D\in\mathbb Z_{\geq 0}$, let
\[
    F_{\le D}(t)
    \colonequals
    \sum_{\substack{\alpha\in\mathbb Z_{\geq 0}^N\\|\alpha|\le D}}
    \frac{c_\alpha}{\alpha!}t^\alpha,
\]
and define
\[
    \|F\|_{\mathcal P^2_{\le D}}^2
    \colonequals
    \|F_{\le D}\|_{\mathcal P^2}^2
    =
    \sum_{\substack{\alpha\in\mathbb Z_{\geq 0}^N\\|\alpha|\le D}}
    \frac{|c_\alpha|^2}{\alpha!}.
\]
We denote by $\mathcal P^2$ the vector space of formal power series
$F$ for which $\|F\|_{\mathcal P^2}<\infty$ (by analogy with the notation for $L^p$ spaces).
The truncated quantity $\|\cdot\|_{\mathcal P^2_{\le D}}$ is a
seminorm on formal power series and a norm on polynomials of total
degree at most $D$.
\end{definition}

The power series viewpoint allows us to exploit the algebraic
structure of these generating functions. The following identity
connects this formalism with the low-degree divergences.

\begin{theorem}[Divergence under Gaussian null model]
\label{thm:gaussian-null}
Let $D\ge1$, and let $\mathbb P$ be a probability measure on $\mathbb R^N$.
Let $\mathcal N=\mathcal N(0,I_N)$. Then
\[
    1+\chi^2_{\le D}(\mathbb P\mid\mathcal N)
    =
    \left\|
        \frac{\phi_{\mathbb P}}{\phi_{\mathcal N}}
    \right\|_{\mathcal P^2_{\le D}}^2,
\]
where the quotient of formal series is also interpreted formally, for terms of total degree up to $D$.
\end{theorem}

\begin{proof}
Let $H_k$ denote the probabilists' Hermite polynomial of degree $k$,
and define
\[
    H_\alpha(y)\colonequals \prod_{i=1}^N H_{\alpha_i}(y_i).
\]
The multivariate Hermite generating function identity is
\[
    \exp\left(\langle t,y\rangle-\frac{\|t\|^2}{2}\right)
    =
    \sum_{\alpha\in\mathbb Z_{\geq 0}^N}
    H_\alpha(y)\frac{t^\alpha}{\alpha!}.
\]
Since
\[
    \phi_{\mathcal N}(t)=\exp\left(\frac{\|t\|^2}{2}\right),
\]
taking expectations coefficientwise through degree $D$ gives
\[
    \left(\frac{\phi_{\mathbb P}}{\phi_{\mathcal N}}\right)_{\le D}(t)
    =
    \sum_{\substack{\alpha\in\mathbb Z_{\geq 0}^N\\|\alpha|\le D}}
    \Ex_{Y \sim \mathbb P}[H_\alpha(Y)]\,\frac{t^\alpha}{\alpha!}.
\]

Moreover,
\[
    \Ex_{Y \sim \mathcal N}[H_\alpha(Y)H_\beta(Y)]
    =
    \alpha!\,\mathbf 1_{\{\alpha=\beta\}}.
\]
Thus every polynomial $f$ of total degree at most $D$ has an
expansion
\[
    f(y)
    =
    \sum_{|\alpha|\le D}
    a_\alpha\frac{H_\alpha(y)}{\sqrt{\alpha!}},
\]
for which
\[
    \Ex_{Y \sim \mathcal N}[f(Y)^2]
    =
    \sum_{|\alpha|\le D}a_\alpha^2
\]
and
\[
    \Ex_{Y \sim \mathbb P}[f(Y)]
    =
    \sum_{|\alpha|\le D}
    a_\alpha
    \frac{\mathbb E_{Y \sim \mathbb P}[H_\alpha(Y)]}{\sqrt{\alpha!}}.
\]
By the definition of $\chi^2_{\le D}$ and
the Cauchy--Schwarz inequality,
\[
    1+\chi^2_{\le D}(\mathbb P\mid\mathcal N)
    =
    \sup_{f\in\mathbb R[Y]_{\le D}\setminus\{0\}}
    \frac{(\mathbb E_{Y \sim \mathbb P} f(Y))^2}
         {\mathbb E_{Y \sim \mathcal N}[f(Y)^2]}
    =
    \sum_{|\alpha|\le D}
    \frac{(\mathbb E_{Y \sim \mathbb P}[H_\alpha(Y)])^2}{\alpha!}.
\]
Comparing this expression with the preceding coefficient identity
proves the result.
\end{proof}

\begin{remark}[Additive Gaussian noise]
The preceding theorem recovers the standard identity for the additive Gaussian noise model~\cite[Theorem~2.6]{kunisky2019notes}.
Suppose that
\[
    Y=X+Z,
    \qquad X\sim\mu,
    \qquad Z\sim\mathcal N,
\]
with $X$ and $Z$ independent. Then $Y\sim\mu*\mathcal N$, and multiplicativity of the formal
MGF gives, through degree $D$,
\[
    \phi_{\mu*\mathcal N}
    =
    \phi_\mu\phi_{\mathcal N}.
\]
Consequently,
\[
    1+\chi^2_{\le D}(\mu*\mathcal N\mid\mathcal N)
    =
    \|\phi_\mu\|_{\mathcal P^2_{\le D}}^2.
\]
\end{remark}

Our analysis of the ratio appearing in Theorem~\ref{thm:gaussian-null}
is based on the elementary identity
\[
    \frac{F(t)}{G(t)}
    =
    \exp\!\left(
        \log F(t)-\log G(t)
    \right).
\]
For $F$ and $G$ MGFs, the logarithms are precisely the
corresponding CGFs. Therefore,
\[
    \frac{\phi_{\mathbb P}(t)}
         {\phi_{\mathbb Q}(t)}
    =
    \exp\!\left(
        \psi_{\mathbb P}(t)-\psi_{\mathbb Q}(t)
    \right).
\]
In particular, when the null distribution is Gaussian,
Theorem~\ref{thm:gaussian-null} may be written equivalently as
\[
    1+\chi^2_{\leq D}(\mathbb P\mid\mathcal N)
    =
    \left\|
        \exp\!\left(
            \psi_{\mathbb P}-\psi_{\mathcal N}
        \right)
    \right\|_{\mathcal P^2_{\leq D}}^2.
\]

This representation is particularly useful because it expresses the
low-degree divergence directly in terms of differences of cumulants.
Consequently, if the low-order cumulants of $\mathbb P$ agree with those
of $\mathcal N$, the corresponding low-degree terms in the exponent
vanish. Our computational lower bounds will follow by controlling the
remaining cumulant coefficients and the combinatorics arising from
exponentiating their generating function.

\subsection{Bounds on truncated generating function norms}
\label{subsec:truncated-generating-function}

We now develop tools for bounding the truncated power series norm
of the exponential of a general formal power series. In our
applications, the exponent will be a difference of
CGFs. We separate the support of its
coefficients from their magnitudes.
The issue of the supports leads to a
combinatorial decomposition problem, while the magnitudes will
be controlled by coarse general moment and cumulant bounds.

To address the supports, we introduce the following set of definitions for a subset of indices $\mathcal B\subseteq\mathbb Z_{\geq 0}^N\setminus\{0\}$.
\begin{definition}
\label{def:decomp-multigraph}
For $\alpha\in\mathbb Z_{\geq 0}^N$, define
\[
    f_{\mathcal B}(\alpha)
    \colonequals
    \max\left(
        \{0\}\cup
        \left\{
            k\ge1:
            \alpha=\beta_1+\cdots+\beta_k,\
            \beta_1,\ldots,\beta_k\in\mathcal B
        \right\}
    \right).
\]
Thus $f_{\mathcal B}(\alpha)=0$ if no such decomposition exists.
\end{definition}

\begin{definition}
\label{def:decomp-class}
For $d\ge1$, define
\[
    \mathcal A_d
    \colonequals
    \left\{
        \alpha\in\mathbb Z_{\geq 0}^N:
        |\alpha|=d,\ f_{\mathcal B}(\alpha)\ge1
    \right\}.
\]
Thus $\mathcal A_d$ consists of the degree-$d$ multi-indices
that admit a decomposition into elements of $\mathcal B$.
\end{definition}

\begin{lemma}
\label{lem:exponential-multiindex-bound}
Let $D\ge1$ be an integer, and let
\[
    F(t)
    =
    \sum_{\beta\in\mathbb Z_{\geq 0}^N\setminus\{0\}}
    \frac{c_\beta}{\beta!}t^\beta
\]
be a formal power series with zero constant term.
Let
$\mathcal B\subseteq\mathbb Z_{\geq 0}^N\setminus\{0\}$
contain the support of its nonzero coefficients:
$c_\beta\ne0 \implies \beta\in\mathcal B$. Suppose that
\[
    |c_\beta|
    \le
    \delta\varepsilon^{|\beta|}(|\beta|!)^K,
    \qquad 1\le|\beta|\le D,
\]
for some $\delta\ge1$, $\varepsilon>0$, and $K\ge1$.
Then
\[
    \|\exp(F)\|_{\mathcal P^2_{\le D}}^2
    \le
    1+
    \sum_{d=1}^{D}
    (4\varepsilon^2)^d(d!)^{2K}
    \sum_{\alpha\in\mathcal A_d}
    \delta^{2f_{\mathcal B}(\alpha)}.
\]
In particular,
\[
    \|\exp(F)\|_{\mathcal P^2_{\le D}}^2
    \le
    1+
    \sum_{d=1}^{D}
    (4\varepsilon^2D^{2K})^d
    \sum_{\alpha\in\mathcal A_d}\delta^{2f_{\mathcal B}(\alpha)}.
\]
\end{lemma}

\begin{proof}
Write
\[
    \exp(F(t))=\sum_{\alpha\in\mathbb Z_{\geq 0}^N}g_\alpha t^\alpha,
    \qquad g_0=1.
\]
If $\alpha\notin\mathcal A_{|\alpha|}$ and $\alpha\ne0$,
then $g_\alpha=0$. Otherwise, fix $\alpha\in\mathcal A_d$
with $d\le D$, and put $f\colonequals f_{\mathcal B}(\alpha)$.
Expanding the exponential and using $\delta\ge1$ gives
\[
\begin{aligned}
    |g_\alpha|
    &\le
    \varepsilon^d
    \sum_{k=1}^{f}\frac{\delta^k}{k!}
    \sum_{\substack{\beta_1+\cdots+\beta_k=\alpha\\
                     \beta_1,\ldots,\beta_k\in\mathcal B}}
    \prod_{j=1}^{k}\frac{(|\beta_j|!)^K}{\beta_j!}\\
    &\le
    \varepsilon^d\delta^f
    \sum_{k=1}^{f}\frac1{k!}
    \sum_{\substack{b_1+\cdots+b_k=d\\b_1,\ldots,b_k\ge1}}
    \left(\prod_{j=1}^{k}(b_j!)^K\right)
    \sum_{\substack{\beta_1+\cdots+\beta_k=\alpha\\
                     |\beta_j|=b_j}}
    \frac1{\beta_1!\cdots\beta_k!}.
\end{aligned}
\]
We grouped by $b_j=|\beta_j|$
and dropped the restrictions $\beta_j\in\mathcal B$. For fixed $b_1,\ldots,b_k$, comparison of the coefficients
of $x^\alpha$ in
\[
    \prod_{j=1}^{k}(x_1+\cdots+x_N)^{b_j}
    =(x_1+\cdots+x_N)^d
\]
yields
\[
    \sum_{\substack{\beta_1+\cdots+\beta_k=\alpha\\
                     |\beta_j|=b_j}}
    \frac1{\beta_1!\cdots\beta_k!}
    =
    \frac{d!}{\alpha!\prod_{j=1}^{k}b_j!}.
\]
Consequently,
\[
    |g_\alpha|
    \le
    \frac{\varepsilon^d\delta^f d!}{\alpha!}
    \sum_{k=1}^{f}\frac1{k!}
    \sum_{\substack{b_1+\cdots+b_k=d\\b_1,\ldots,b_k\ge1}}
    \prod_{j=1}^{k}(b_j!)^{K-1}.
\]
Since $K\ge1$ and $\prod_jb_j!\le d!$, and there are
$\binom{d-1}{k-1}$ positive ordered compositions of $d$
into $k$ parts, we obtain
\[
\begin{aligned}
    |g_\alpha|
    &\le
    \frac{\varepsilon^d\delta^f(d!)^K}{\alpha!}
    \sum_{k=1}^{f}\frac1{k!}\binom{d-1}{k-1}\\
    &\le
    \frac{2^{d-1}\varepsilon^d(d!)^K}{\alpha!}
    \delta^{f_{\mathcal B}(\alpha)}.
\end{aligned}
\]
Therefore,
\[
\begin{aligned}
    \|\exp(F)\|_{\mathcal P^2_{\le D}}^2
    &=1+\sum_{1\le|\alpha|\le D}\alpha!|g_\alpha|^2\\
    &\le
    1+\sum_{d=1}^{D}
    4^{d-1}\varepsilon^{2d}(d!)^{2K}
    \sum_{\alpha\in\mathcal A_d}
    \frac{\delta^{2f_{\mathcal B}(\alpha)}}{\alpha!} \\
    &\le 1+\sum_{d=1}^{D}
    (2\varepsilon)^{2d}(d!)^{2K}
    \sum_{\alpha\in\mathcal A_d}
    \delta^{2f_{\mathcal B}(\alpha)}
\end{aligned}
\]
The second result follows from $d!\le D^d$ for $d\le D$.
\end{proof}

\subsection{Toolkit for cumulants}

The above reasoning addresses the combinatorial side of the terms arising when exponentiating a possibly sparse generating function.
We now give some facts and estimates concerning cumulants that will let us bound the required terms in our particular application.
We write $\kappa(U_1,\ldots,U_d)$ for the joint cumulant of a finite family of real random variables.
For $X\sim\mu$ and $U_j=X_{i_j}$, this agrees with $\kappa_\alpha(\mu)$, where $\alpha_i$ counts the occurrences of $i$ among $i_1,\ldots,i_d$.
The following identities also define joint cumulants when only the indicated mixed moments exist.

\begin{proposition}[Moment--cumulant identities]
\label{prop:moment-cumulant-identities}
Let $d\ge1$ and let $U_1,\ldots,U_d$ be jointly distributed real random variables such that
$\mathbb{E}[\prod_{i\in S}|U_i|]<\infty$ for every nonempty $S\subseteq[d]$.
Write $\Part(S)$ for the set of partitions of $S$ into nonempty blocks.
Then
\begin{align*}
    \kappa(U_1,\ldots,U_d)
    &=\sum_{\pi\in\Part([d])}(-1)^{|\pi|-1}(|\pi|-1)!
      \prod_{B\in\pi}\mathbb{E}\left[\prod_{i\in B}U_i\right],\\
    \mathbb{E}\left[\prod_{i=1}^d U_i\right]
    &=\sum_{\pi\in\Part([d])}\prod_{B\in\pi}\kappa(U_i:i\in B).
\end{align*}
In particular, the first and second cumulants are the mean and covariance, and the third joint cumulant of centered variables is their product moment.
Joint cumulants are symmetric and multilinear in their inputs.
\end{proposition}

\begin{proof}
It suffices to use the moment polynomial
\[
    \Phi(t)\colonequals
    \sum_{S\subseteq[d]}\mathbb{E}\left[\prod_{i\in S}U_i\right]\prod_{i\in S}t_i,
\]
with empty products equal to one.
For nonempty $S\subseteq[d]$, the coefficient of $\prod_{i\in S}t_i$ in the formal logarithm $\log\Phi(t)$ is $\kappa(U_i:i\in S)$.
Expanding $\log\Phi$ gives the first identity: a partition into $k$ blocks occurs in $k!$ orders, with coefficient $(-1)^{k-1}/k$.
Expanding $\Phi=\exp(\log\Phi)$ gives the second identity.
The remaining claims follow from the first formula.
\end{proof}

\begin{proposition}[Vanishing across independent groups]
\label{prop:cumulant-independence}
Under the moment assumptions of Proposition~\ref{prop:moment-cumulant-identities}, suppose $[d]=A\sqcup B$ with $A,B$ nonempty and the random vectors $(U_i)_{i\in A}$ and $(U_i)_{i\in B}$ independent.
Then
\[
    \kappa(U_1,\ldots,U_d)=0.
\]
\end{proposition}

\begin{proof}
Independence gives $\Phi=\Phi_A\Phi_B$ for the moment polynomial in the preceding proof, where each factor uses only the variables in its group.
Thus $\log\Phi=\log\Phi_A+\log\Phi_B$ has no term involving all $d$ variables.
\end{proof}

We next give a general-purpose bound on cumulants in terms of moments under the following very mild growth property.
\begin{definition}[$(K,C)$-bounded moments]
\label{def:bounded-moments}
Let $K\ge1$ and $C>0$. A probability measure $\mu$
on $\mathbb R^N$ with moments of all orders has
\emph{$(K,C)$-bounded moments} if, for every $\alpha\in\mathbb Z_{\geq 0}^N$,
\[
    |m_\alpha(\mu)|
    =
    \left|\mathbb E_{X\sim\mu}X^\alpha\right|
    \le C^{|\alpha|}(|\alpha|!)^K.
\]
\end{definition}
\noindent
The following then bounds cumulants in the case of
$(K,C)$-bounded moments.

\begin{lemma}[A priori cumulant bound]
\label{lem:apriori-cumulant-bound}
Let $\mu$ be a probability measure on $\mathbb R^N$ with
$(K,C)$-bounded moments. Then, for every nonzero
$\alpha\in\mathbb Z_{\geq 0}^N$,
\[
    |\kappa_\alpha(\mu)|
    \le (2C)^{|\alpha|}(|\alpha|!)^K.
\]
\end{lemma}

\begin{proof}
Set $d\colonequals |\alpha|$.
Choose indices $i_1,\ldots,i_d\in[N]$ such that each
$i\in[N]$ occurs exactly $\alpha_i$ times.
For each nonempty block $B\subseteq[d]$, let
$\beta(B)\in\mathbb Z_{\geq 0}^N$ record the multiplicities of
the indices $(i_j)_{j\in B}$. Then $|\beta(B)|=|B|$, so
\[
    \left|
        \mathbb E_{X\sim\mu}\prod_{j\in B}X_{i_j}
    \right|
    =
    |m_{\beta(B)}(\mu)|
    \le
    C^{|B|}(|B|!)^K.
\]

Proposition~\ref{prop:moment-cumulant-identities} and the preceding bound give
\[
    |\kappa_\alpha(\mu)|
    \le
    C^d
    \sum_{\pi\in\operatorname{Part}([d])}
    (|\pi|-1)!
    \prod_{B\in\pi}(|B|!)^K.
\]
For a fixed number of blocks $k$, ordering the blocks and
grouping according to their sizes yields
\[
    \sum_{\substack{
        \pi\in\operatorname{Part}([d])\\
        |\pi|=k}}
    (k-1)!\prod_{B\in\pi}(|B|!)^K =
    \frac{d!}{k}
    \sum_{\substack{
        d_1+\cdots+d_k=d\\
        d_1,\ldots,d_k\ge1}}
    \prod_{j=1}^k(d_j!)^{K-1}.
\]
Since $K\ge1$ and
\[
    \prod_{j=1}^k d_j!\le d!,
\]
each product in the last sum is at most $(d!)^{K-1}$.
There are $\binom{d-1}{k-1}$ compositions of $d$ into
$k$ positive parts. Therefore,
\[
    \sum_{\substack{
        \pi\in\operatorname{Part}([d])\\
        |\pi|=k}}
    (k-1)!\prod_{B\in\pi}(|B|!)^K
    \le
    \frac{(d!)^K}{k}\binom{d-1}{k-1}.
\]
Summing over $k$ gives
\begin{align*}
    |\kappa_\alpha(\mu)|
    &\le
    C^d(d!)^K
    \sum_{k=1}^d\frac{1}{k}\binom{d-1}{k-1}
    \\
    &=
    C^d(d!)^K\frac{1}{d}
    \sum_{k=1}^d\binom{d}{k}
    \\
    &=
    C^d(d!)^K\frac{2^d-1}{d}
    \\
    &\le
    (2C)^d(d!)^K. \qedhere
\end{align*}
\end{proof}

\begin{proposition}[Cumulants of CLT-normalized convolutions]
\label{prop:normalized-convolution-cumulants}
Let $\mu$ be a probability measure on $\mathbb R^N$ with
moments of every order, and let $r\ge1$ be an integer.
Suppose $\mu$ has $(K,C)$-bounded moments with $K\ge1$.
Then, for every nonzero multi-index $\alpha \in \mathbb Z_{\geq 0}^N$,
\[
    |\kappa_\alpha(r^{-1/2}\cdot \mu^{*r})|
    \le
    (2C)^{|\alpha|} r^{1 - |\alpha| /2}(|\alpha|!)^K.
\]
\end{proposition}

\begin{proof}
The bound follows by applying Lemma~\ref{lem:apriori-cumulant-bound} to the cumulant on the right-hand side of Proposition~\ref{prop:scaling-cumulants}.
\end{proof}

\begin{remark}
It follows immediately by the same arguments that, for conclusions involving only cumulants through degree $D$,
it suffices to assume finite moments and the stated moment
bounds through degree $D$.
\end{remark}

The following result combines the ingredients established above to reduce the control of the low-degree
divergence for a CLT testing problem to the combinatorics of the support of the
non-Gaussian cumulants.

\begin{lemma}[Low-degree bound from cumulant support]
\label{lem:clt-cumulant-support-bound}
Let $\mu$ be a centered and isotropic probability measure on
$\mathbb R^N$ with $(K,C)$-bounded moments, where $K\ge1$.
Let
$\varnothing\ne\mathcal B \subseteq\{\beta\in\mathbb Z_{\geq 0}^N:|\beta|\ge3\}$
satisfy $\kappa_\beta(\mu)=0$ whenever $|\beta|\ge3$
and $\beta\notin\mathcal B$, and set
$q\colonequals \min_{\beta\in\mathcal B}|\beta|$. For all
integers $r\ge1$ and $D\ge q$,
\[
    \chi^2_{\le D}(r^{-1/2}\cdot \mu^{*r}\mid\mathcal N)
    \le
    \sum_{d=q}^D
    \left(\frac{16C^2D^{2K}}{r}\right)^d
    \sum_{\alpha \in \mathcal{A}_d}
    r^{2f_{\mathcal B}(\alpha)}.
\]
\end{lemma}

\begin{proof}
Set
\[
    F(t)\colonequals
    \psi_{r^{-1/2}\cdot \mu^{*r}}(t)-\psi_{\mathcal N}(t).
\]
Since $\mu$ is centered and isotropic, the terms of degrees
one and two cancel. Thus
\[
    F(t)=
    \sum_{\beta\in\mathcal B}
    \frac{c_\beta}{\beta!}t^\beta,
    \qquad
    c_\beta=r^{1-|\beta|/2}\kappa_\beta(\mu).
\]
By Proposition~\ref{prop:normalized-convolution-cumulants},
\[
    |c_\beta|
    \le
    r\left(\frac{2C}{\sqrt r}\right)^{|\beta|}
    (|\beta|!)^K.
\]
We apply Lemma~\ref{lem:exponential-multiindex-bound} with
$\delta=r$ and $\varepsilon=2C/\sqrt r$.
Using Theorem~\ref{thm:gaussian-null} gives
\[
\begin{aligned}
    \chi^2_{\le D}
    \left(r^{-1/2}\cdot \mu^{*r}\mid\mathcal N\right)
    &=
    \|\exp(F)\|_{\mathcal P^2_{\le D}}^2-1
    \\
    &\le
    \sum_{d=q}^D\sum_{\alpha\in\mathcal A_d}
    \left(\frac{16C^2D^{2K}}{r}\right)^{|\alpha|}
    r^{2f_{\mathcal B}(\alpha)},
\end{aligned}
\]
where all terms with $0<|\alpha|<q$ vanish because no such
$\alpha$ admits a decomposition into elements of $\mathcal B$.
\end{proof}

Finally, we give a last tool that allows us to convert moment bounds on the distributions $\nu$ of vector entries to those on the summands $X$ in tensor CLT models.
\begin{lemma}[Moment bounds for tensor summands]
\label{lem:tensor-moment-bound}
Let $\mu_n$ be the tensor summand distribution in the
Wigner--Wishart tensor testing problem of Definition~\ref{def:wigner-wishart}, and suppose that
$\nu$ has $(K,C)$-bounded moments with $K\ge1$.
Then $\mu_n$ has $(pK,C^p)$-bounded moments, uniformly in $n$.
\end{lemma}

\begin{proof}
Fix $\alpha\in\mathbb Z_{\geq 0}^{\binom np}$, write $d\colonequals |\alpha|$, and, for each $i\in[n]$, define
\[
    d_i\colonequals \sum_{\substack{e\in\binom{[n]}p: \ i\in e}}\alpha_e.
\]
By independence of the coordinates of the underlying vector,
\[
    m_\alpha(\mu_n)
    =
    \prod_{i=1}^n m_{d_i}(\nu).
\]
The moment assumption on $\nu$ therefore gives
\[
    |m_\alpha(\mu_n)|
    \le
    C^{\sum_i d_i}\prod_{i=1}^n(d_i!)^K
    =
    C^{pd}\prod_{i=1}^n(d_i!)^K.
\]
If $d\ge1$, then $0\le d_i\le d$ and $\sum_i d_i=pd$.
Since $(m!)^{1/m}$ is nondecreasing for integers $m\ge1$,
\[
    \prod_{i=1}^n d_i!
    \le
    \prod_{i=1}^n(d!)^{d_i/d}
    =
    (d!)^p.
\]
Hence
\[
    |m_\alpha(\mu_n)|
    \le
    (C^p)^d(d!)^{pK}. \qedhere
\]
\end{proof}

\section{Proofs of main results}

\subsection{Statistical lower bounds from total variation distances}

We first develop the estimates for the undiluted model and prove Theorem~\ref{thm:full-tensor-tv}.
We then adapt the argument to prove the diluted bound in Theorem~\ref{thm:diluted-tensor-tv}.
Throughout this section, we work in the model of Definition~\ref{def:wigner-wishart}, $p$ is fixed, and $\nu=\mathcal N(0,1)$ unless explicitly stated otherwise.
The notation $a\lesssim b$ means $a\le Cb$ for a positive constant $C$.
Implicit constants in $\lesssim$, as well as explicit constants $C,c>0$, may change from line to line and may depend on $p,\nu$ and any parameters indicated by subscripts, but not on $n$ or $r$.
In Section~\ref{sec:pf:thm:diluted-tensor-tv}, $c$ denotes the dilution parameter, and the implicit constants are also independent of $c$.
We assume $n\ge p$, since otherwise there are no tensor coordinates.
Write $\mathcal I=\binom{[n]}p$ and $N=|\mathcal I|$.
Let $X=X(v)\sim\mu_n$ denote one summand, and for $r\ge1$ set
\[
    Y = Y_r=r^{-1/2}\sum_{j=1}^r X(v^{(j)}),
\]
where $v,v^{(1)},\ldots,v^{(r)}$ are independent with law $\nu^{\otimes n}$.
Thus $\mathbb P_n=\Law(Y_r)$ and $\mathbb Q_n=\mathcal N = \mathcal N(0,I_N)$.
We write $\|T\|_{\mathrm{F}}^2$ for the squared Frobenius norm of a tensor $T$, the sum of the squares of all its entries, and $\|T\|_{L^2(\mathrm{F})}^2=\mathbb{E}[\|T\|_{\mathrm{F}}^2]$ for a random tensor.
Indices $A,B,C,A_j$ below range over $\mathcal I$.

\subsubsection{Overview}
\label{sec:overview}

We give an overview of the main proof ideas for the undiluted model, as this is the most technical argument we present and also the one whose discovery involved the most AI assistance (as mentioned in Section~\ref{sec:ai}).
Recall that we wish to bound the TV distance $d_{\mathrm{TV}}(\mathbb P_n, \mathcal N)$.
The general technique is \emph{Stein's method}, which is based on characterizing $\mathcal N$ as the unique law of a random vector $Y \in \mathbb{R}^N$ such that Stein's identity
\[ \mathcal{E}(f) = \mathcal{E}_r(f) \colonequals \mathbb{E}[\langle Y, f(Y) \rangle - \operatorname{div} f(Y)] = 0 \]
holds exactly for all smooth compactly supported $f: \mathbb{R}^N \to \mathbb{R}^N$ (this is the same as the ``Gaussian integration by parts'' identity).
The general premise of the method is that, even if Stein's identity does not hold exactly, if the residual $\mathcal{E}(f)$ is small for all such test functions $f$ then $\Law(Y) \approx \mathcal{N}$.
The following gives a general concrete version of this, following from standard inequalities relating Fisher information, relative entropy, and distributional distances recalled in~\cite[Section~1]{LNP-2015-SteinLogSobolevTransport}.
\begin{lemma}[Distances from the Gaussian Stein identity]
\label{lem:stein-error-distance}
Let $Y\sim\mathbb{P}$ have a density on $\mathbb R^N$ and finite second moment.
Assume that $\mathbb{P}$ has a weak score function $\rho \in L^2(\mathbb{P})$, meaning that $\mathbb E[\nabla g(Y)]=-\mathbb E[g(Y)\rho(Y)]$ for every $g\in C_c^\infty(\mathbb R^N)$.
For $f\in C_c^\infty(\mathbb R^N;\mathbb R^N)$, set
\begin{align*}
    \mathcal E(f)&\colonequals
    \mathbb E[\langle Y,f(Y)\rangle-\operatorname{div}f(Y)],
    \\
    \varepsilon&\colonequals
    \sup_{\substack{f\in C_c^\infty(\mathbb R^N;\mathbb R^N)\\
                    \|f(Y)\|_{L^2}^2 = \mathbb E[\|f(Y)\|_2^2]\le1}}
    |\mathcal E(f)|.
\end{align*}
Then
\[
    \varepsilon^2=I(\mathbb{P}\mid\mathcal N)
    \colonequals\mathbb E[\|Y+\rho(Y)\|_2^2],
\]
and
\[
    d_{\mathrm{TV}}(\mathbb{P},\mathcal N)\le\frac{\varepsilon}{2},
    \qquad
    W_2(\mathbb{P},\mathcal N)\le\varepsilon,
    \qquad
    \operatorname{Ent}(\mathbb{P}\mid\mathcal N)\le\frac{\varepsilon^2}{2}.
\]
Here $I$ is the relative Fisher information, $W_2$ is the quadratic Wasserstein distance for the Euclidean norm, and
$\operatorname{Ent}(\mathbb{P}\mid\mathcal N)\colonequals\int\log(d\mathbb{P}/d\mathcal N)\,d\mathbb{P}$ is relative entropy.
\end{lemma}
\noindent
We note that the approach of \cite{Mikulincer-2020-CLTWishartTensors} is similar in that it uses Stein's method and analyzes the error functional $\mathcal{E}(f)$, but it controls $|\mathcal{E}(f)|$ relative to a norm of \emph{derivatives} of $f$, namely $\|Df(Y)\|_{L^2}$, leading to control only of the Wasserstein distance.
Since we control $|\mathcal{E}(f)|$ relative to the norm of $f$ itself, we obtain bounds on the relative Fisher information and therefore the TV distance as well.
This route via the Fisher information was already suggested as a technically challenging alternative in Section~2.1.1 of \cite{Mikulincer-2020-CLTWishartTensors}.

Our approach is as follows.
First, recall that our $Y$ is itself a polynomial of underlying Gaussian random variables $v^{(1)}, \dots, v^{(r)}$.
This simplifies our task, because we may use Gaussian integration by parts exactly with respect to \emph{these} random variables to obtain the approximate version of Gaussian integration by parts that $Y$ satisfies.
The first step of our plan is to use this approach to produce an expansion of $\mathcal{E}(f)$ for a smooth compactly supported $f$ in terms of the joint cumulants of $Y$ and a remainder term (Lemma~\ref{lem:tv-cumulant-expansion}).

However, we do not directly bound $\mathcal{E}(f)$ using this expansion, but rather apply the expansion on only a subset of summands of $Y$.
Namely, suppose we can write $r = Ks$ (otherwise we take our subsets to be as close to equally sized as possible) and write $W_1, \dots, W_{K}$ for the suitably normalized sums of $s$ of the terms in $Y$, so that $\Law(W_i) = \Law(Y_s)$.
Then, we have $Y = \frac{1}{\sqrt{K}}\sum_{i = 1}^{K} W_i$, and the $W_i$ are independent.
Writing $\mathbb{E}_{\sim i}$ for the expectation with respect to all $W_j$ for $j \neq i$, we then rewrite
\begin{equation}
    \mathcal{E}_r(f) = \frac{1}{\sqrt{K}}\sum_{i = 1}^K \mathbb{E}_{\sim i}[\mathcal{E}_s(f_{(i)})], \label{eq:block-expansion}
\end{equation}
where $f_{(i)}(w) \colonequals f(\frac{1}{\sqrt{K}}\,w + \frac{1}{\sqrt{K}}\sum_{j\ne i}W_j)$ after conditioning on the $(W_j)_{j \neq i}$.

When we apply our cumulant expansion to the $\mathcal E_s(f_{(i)})$ terms in \eqref{eq:block-expansion}, we end up with terms depending on $W_i$ multiplied by derivatives of $f(Y)$, then averaged over all other $W_j$.
We integrate by parts with respect to these other $W_j$, the outcome of which may be formulated in terms of replacing the derivative terms by \emph{score} terms,
\[ \rho(w) = \rho_s(w) = \nabla \log q_s(w) \in \mathbb{R}^N, \]
where $q_s(w)$ is the density of each $W_j$ (equivalently, of $Y_s$), and the gradient is understood weakly as in Lemma~\ref{lem:stein-error-distance}.
Thus after this manipulation we replace derivatives of $f(Y)$ by products of components of $\rho_s(W_j)$ for various $W_j$.
Lemma~\ref{lem:tv-block-representation} writes the final result in a form that can be bounded by the Cauchy-Schwarz inequality as
\[ |\mathcal{E}_r(f)| = \left|\mathbb{E}[\langle H(W_1, \dots, W_K), f(Y) \rangle]\right| \leq \|H(W_1, \dots, W_K)\|_{L^2} \|f(Y)\|_{L^2}, \]
where $H = H(W_1, \dots, W_K)$ is a complicated random vector whose entries are various expressions in the scores $\rho_s(W_j)$ and the entries of the $W_j$.
Crucially, $H$ does not depend on $f$, so in the context of Lemma~\ref{lem:stein-error-distance} we have $\epsilon \leq \|H\|_{L^2}$ and it remains to control this single norm.
Also importantly, terms in entries of $H$ are multilinear in the scores $\rho_s(W_j)$ for different $j$, and thus these parts of the terms in $H$ are products of independent such factors; to ensure this, we integrate by parts against \emph{different} $W_j$ when replacing each derivative term above.
This requires that the number of blocks $K$ exceed the highest derivative order $k$ appearing in the cumulant expansion of $\mathcal{E}_s(f)$ we use; we simply take $K=k+1$.

Finally, we control $\|H\|_{L^2}$ by controlling the score functions.
In particular, using the above observations of independence and multilinearity in the scores, it turns out to be enough to control the second moments of the score vectors with an estimate of the form
\[ \mathbb E[\rho_s(W_j)\rho_s(W_j)^\top] \preceq C_pI_N. \]
This is the key analytic estimate of the proof, given in Lemma~\ref{lem:tv-block-score} provided $s \geq C_pn^p$ for sufficiently large $C_p$.
Repeatedly applying this to bound the scores and combining with estimates on the cumulants and remainder terms from the initial expansion, we may then fairly directly bound $\|H\|_{L^2}$, and the main result then follows from Lemma~\ref{lem:stein-error-distance}.

\subsubsection{Cumulant expansion of Stein identity error}
\label{sec:tv-gaussian-kernels}

Following the above plan, we first expand $\mathcal E_r(f)$ in cumulants and a remainder by repeated Gaussian integration by parts.
Related higher-order Stein identities and their applications to CLTs are developed in~\cite{Fathi-2021-HigherOrderSteinKernels}.

For $m\ge2$, write $\kappa_m(X)$ for the array with entries
\[
    \kappa_m(X)_{A_1,\ldots,A_m}
    \colonequals\kappa(X_{A_1},\ldots,X_{A_m}).
\]
To distinguish polynomials from their random values, write
\[
    x_A(v)\colonequals\prod_{i\in A}v_i,
    \qquad X_A=x_A(v).
\]
All derivatives in the construction below are with respect to the underlying coordinates $v_i$.
In particular,
\[
    \partial_i x_A(v)=\bm 1\{i\in A\}\prod_{\ell\in A\setminus\{i\}}v_\ell.
\]
Let $\mathcal L=\sum_{i=1}^n(\partial_i^2-v_i\partial_i)$ be the Ornstein--Uhlenbeck generator acting on these polynomials \cite[Section~2.7]{BGL-2014-MarkovDiffusionOperators}.
Gaussian integration by parts gives
\[
    -\mathbb{E}[f(v)(\mathcal Lg)(v)]
    =\sum_{i=1}^n\mathbb{E}[\partial_i f(v)\partial_i g(v)].
\]
For a centered polynomial $h$, define $\mathcal Rh=(-\mathcal L)^{-1}h$ by its Hermite expansion:
\[
    h(v)=\sum_{|\alpha|>0}h_\alpha H_\alpha(v),
    \qquad
    (\mathcal Rh)(v)=\sum_{|\alpha|>0}\frac{h_\alpha}{|\alpha|}H_\alpha(v),
    \qquad H_\alpha(v)=\prod_i H_{\alpha_i}(v_i).
\]
Here $H_j$ is the probabilists' Hermite polynomial, so $\mathcal LH_\alpha=-|\alpha|H_\alpha$.
The inverse $\mathcal R$ preserves polynomial degree and dependence on a given set of vector coordinates, and
\begin{equation}
\label{eq:tv-resolvent-ibp}
    \mathbb{E}[h(v)g(v)]=\sum_{i=1}^n\mathbb{E}[\partial_i(\mathcal Rh)(v)\partial_i g(v)].
\end{equation}
In particular, $\mathcal Lx_A=-px_A$ and $\mathcal Rx_A=x_A/p$.

We construct tensors $T_k(v)$ whose entries are polynomials in $v$, and write $M_k\colonequals\mathbb{E}[T_k(v)]$ for their entrywise means.
The first tensor is
\begin{equation}
\label{eq:tv-kernel-first}
    (T_1(v))_{A,B}
    =\frac1p\sum_{i\in A\cap B}
      \left(\prod_{\ell\in A\setminus\{i\}}v_\ell\right)
      \left(\prod_{\ell\in B\setminus\{i\}}v_\ell\right)
      -\delta_{A,B}.
\end{equation}
This is $p^{-1}\sum_i\partial_i x_A(v)\partial_i x_B(v)-\delta_{A,B}$, so \eqref{eq:tv-resolvent-ibp} and $\Cov(X)=I_N$ give $M_1=0$.

For the next tensor, fix $A_0,\ldots,A_{k+1}$.
For each $j\in\{1,\ldots,k+1\}$, omit $A_j$ from the last $k+1$ indices, subtract the mean of the resulting entry of $T_k$, and apply $\mathcal R$ to that scalar polynomial:
\[
    h_j\colonequals
    \mathcal R\left[(T_k-M_k)_{A_0,\ldots,\widehat{A_j},\ldots,A_{k+1}}\right].
\]
Here the hat means that the indicated index is omitted.
Define
\begin{equation}
\label{eq:tv-kernel-recursion}
    (T_{k+1}(v))_{A_0,\ldots,A_{k+1}}
    =\frac1{k+1}\sum_{j=1}^{k+1}\sum_{i\in A_j}
      \partial_i h_j(v)\prod_{\ell\in A_j\setminus\{i\}}v_\ell.
\end{equation}
The inner sum is the ordinary gradient inner product $\langle\nabla h_j(v),\nabla x_{A_j}(v)\rangle$.
Since $T_k$ is symmetric in its last $k$ indices, averaging over the $k+1$ choices of $j$ makes $T_{k+1}$ symmetric in its last $k+1$ indices, with $A_0$ fixed.
Mixed derivatives are symmetric in the same indices, so averaging preserves the sums against these derivatives used below.

The cumulant identity in the next proposition is a multivariate counterpart of the recursive formula in~\cite[Theorem~4.3]{NP-2010-CumulantsWienerSpace}.

\begin{proposition}[Fluctuation bounds and means]
\label{prop:tv-kernels}
For each fixed $k\ge1$, we have
\begin{equation}
\label{eq:tv-kernel-bound}
    \|T_k(v)-M_k\|_{L^2(\mathrm{F})}
    \lesssim_{p,k,\nu} n^{((k+1)p-k)/2}.
\end{equation}
For $k\ge2$, we also have
\begin{equation}
    M_k = \frac{1}{k!}\kappa_{k+1}(X).
\end{equation}
\end{proposition}
\begin{proof}
Each new derivative contraction requires the appended $p$-set to meet the union of the preceding sets.
Thus a nonzero entry of $T_k$ has
$|A_0\cup\cdots\cup A_k|\le(k+1)p-k$.
Each entry is a polynomial of bounded degree in a bounded number of Gaussian coordinates.
The Hermite formula for $\mathcal R$ gives a uniform $L^2$ bound for each overlap pattern, and only finitely many patterns occur for fixed $p,k$.
Counting these patterns proves \eqref{eq:tv-kernel-bound}.

To identify the mean, fix $A_0$ and define the polynomial $F(v)=\sum_A u_Ax_A(v)$.
Set $G_1(v)=\sum_i\partial_i(\mathcal Rx_{A_0})(v)\partial_i F(v)$.
For $j\ge1$, set
\[
    G_{j+1}(v)=\sum_i
      \partial_i\bigl[\mathcal R(G_j-\mathbb E[G_j(v)])\bigr](v)
      \partial_iF(v).
\]
For $j\ge2$, $G_j$ is the contraction of $T_j$ with $u^{\otimes j}$ in its last $j$ indices, with first index $A_0$.
Repeated use of \eqref{eq:tv-resolvent-ibp} gives
\[
    \mathbb{E}[X_{A_0}F(v)^k]
    =\sum_{j=1}^k\frac{k!}{(k-j)!}\mathbb{E}[G_j(v)]\mathbb{E}[F(v)^{k-j}].
\]
The moment--cumulant identity in Proposition~\ref{prop:moment-cumulant-identities}, with $X_{A_0}$ distinguished, expresses the same moment as
\[
    \sum_{j=1}^k\binom{k}{j}
    \kappa(X_{A_0},\underbrace{F(v),\ldots,F(v)}_{j\text{ copies}})\mathbb{E}[F(v)^{k-j}].
\]
Induction gives $\mathbb{E}[G_j(v)]=\kappa(X_{A_0},F(v),\ldots,F(v))/j!$.
Comparing coefficients of the $u_A$ proves the claimed identity for the tensors, which are symmetric in their last $j$ indices.
\end{proof}

For a smooth compactly supported $f:\mathbb R^N\to\mathbb R^N$, define the derivative tensor
\[
    (D^mf)_{A_0,\ldots,A_m}=\partial_{A_1}\cdots\partial_{A_m}f_{A_0}.
\]
Here $\partial_A$ denotes differentiation with respect to the observation coordinate $y_A$.
We use the entrywise inner product for array contractions.

\begin{lemma}[Cumulant expansion of Stein identity error]
\label{lem:tv-cumulant-expansion}
In the above setting, for every fixed $k\ge2$, every $r\ge1$, and every $f\in C_c^\infty(\mathbb R^N;\mathbb R^N)$,
\begin{equation}
\label{eq:tv-block-identity}
    \mathcal E_r(f)
    =\sum_{m=2}^k r^{-(m-1)/2}
      \langle M_m,\mathbb{E}[D^mf(Y_r)]\rangle
      +\mathbb{E}[\langle R_{k,r}(Y_r),D^kf(Y_r)\rangle],
\end{equation}
where $M_m=\kappa_{m+1}(X)/m!$ and
\begin{equation}
\label{eq:tv-remainder-definition}
    R_{k,r}(Y_r)
    \colonequals r^{-(k+1)/2}\mathbb{E}\left[
      \sum_{j=1}^r(T_k(v^{(j)})-M_k)\,\middle|\,Y_r\right].
\end{equation}
This remainder satisfies
\begin{equation}
\label{eq:tv-block-remainder}
    \|R_{k,r}(Y_r)\|_{L^2(\mathrm{F})}
    \lesssim_{p,k}\frac{n^{((k+1)p-k)/2}}{r^{k/2}}.
\end{equation}
\end{lemma}
\begin{proof}
Applying \eqref{eq:tv-resolvent-ibp} to each underlying vector gives
\[
    \mathcal E_r(f)
    =\frac1r\sum_{j=1}^r
      \mathbb{E}[\langle T_1(v^{(j)}),Df(Y_r)\rangle].
\]
At each further integration, we separate the mean of the coefficient tensor from its centered part and apply \eqref{eq:tv-resolvent-ibp} to the latter.
The chain rule introduces a factor $r^{-1/2}$ and the recursion~\eqref{eq:tv-kernel-recursion} gives the next tensor.
After reaching $T_k$, we condition the remaining centered coefficient on $Y_r$ to obtain \eqref{eq:tv-block-identity} and \eqref{eq:tv-remainder-definition}.
Proposition~\ref{prop:tv-kernels} identifies the means as cumulants.
Using a conditional Jensen inequality and independence gives
\[
    \|R_{k,r}(Y_r)\|_{L^2(\mathrm{F})}
    \le r^{-k/2}\|T_k(v)-M_k\|_{L^2(\mathrm{F})},
\]
so the fluctuation bound in the same proposition proves \eqref{eq:tv-block-remainder}.
\end{proof}

\subsubsection{Transferring derivatives to independent scores}
\label{sec:tv-block-transfer}

We next construct the vector $H$ from Section~\ref{sec:overview}.
We assume a bound on the score functions involved here and prove it in the following section.
For a random vector $W$ admitting a weak score $\rho_W$ in the sense of Lemma~\ref{lem:stein-error-distance}, we have
$\mathbb{E}[\nabla h(W)]=-\mathbb{E}[h(W)\rho_W(W)]$ for smooth compactly supported $h$.

\begin{lemma}[Block representation of Stein identity error]
\label{lem:tv-block-representation}
Fix $k\ge2$ and $r\ge k+1$.
Partition the $r$ summands in $Y = Y_r$ into $K=k+1$ groups of positive sizes $r_i$ differing by at most one.
Let $W_i$ be their normalized sums, and set $c_i=\sqrt{r_i/r}$, so that
\[
    W_i\stackrel{d}{=}Y_{r_i},\qquad
    Y_r=\sum_{i=1}^K c_iW_i,\qquad
    \sum_{i=1}^K c_i^2=1.
\]
Suppose that the expansion~\eqref{eq:tv-block-identity} and remainder bound~\eqref{eq:tv-block-remainder} hold at every block size $r_i$.
Suppose also that each $W_i$ has a density and a weak score $\rho_i$ satisfying
\[
    \mathbb{E}[\rho_i(W_i)\rho_i(W_i)^\top]\preceq L I_N
\]
for a constant $L$.
Then the vector $H=H(W_1,\ldots,W_K)$ defined in~\eqref{eq:tv-H-definition} below is square-integrable, independent of the choice of $f$, and satisfies
\begin{equation}
\label{eq:tv-H-representation}
    \mathcal E_r(f)=\mathbb{E}[\langle H,f(Y_r)\rangle]
    \qquad\text{for every }f\in C_c^\infty(\mathbb R^N;\mathbb R^N),
\end{equation}
and
\begin{equation}
\label{eq:tv-H-bound}
    \|H\|_{L^2}
    \lesssim_{p,k,L}\left(
       \sum_{m=2}^k\frac{\|M_m\|_{\mathrm F}}{r^{(m-1)/2}}
       +\frac{n^{((k+1)p-k)/2}}{r^{k/2}}\right).
\end{equation}
Each term in $H$ is a contraction of a deterministic cumulant tensor or a remainder depending only on $W_i$ with scores from distinct blocks other than $i$.
\end{lemma}
\begin{proof}
We condition on all blocks except $W_i$ and apply \eqref{eq:tv-block-identity} to
\[
    f_{(i)}(w)=f\left(c_iw+\sum_{j\ne i}c_jW_j\right).
\]
Multiplying by $c_i$ and summing over $i$ gives
\begin{align*}
    \mathcal E_r(f)
    &=\sum_{i=1}^K\sum_{m=2}^k
      c_i^{m+1}r_i^{-(m-1)/2}
      \mathbb{E}[\langle M_m,D^mf(Y_r)\rangle]\\
    &\hspace{1cm}+\sum_{i=1}^K c_i^{k+1}
      \mathbb{E}[\langle R_{k,r_i}(W_i),D^kf(Y_r)\rangle].
\end{align*}
Here the divergence terms combine because $\sum_i c_i^2=1$.

For each $i$ and $2\le m\le k$, choose distinct indices $j_1,\ldots,j_m$ other than $i$, and put
\[
    U_{m,i}=\rho_{j_1}(W_{j_1})\otimes\cdots\otimes\rho_{j_m}(W_{j_m}),
    \qquad
    a_{m,i}=(-1)^m c_i^{m+1}\prod_{h=1}^m c_{j_h}^{-1}.
\]
Integration by parts in block $j_h$ replaces $\partial_A$ by $-c_{j_h}^{-1}\rho_{j_h,A}(W_{j_h})$.
Thus \eqref{eq:tv-H-representation} holds with
\begin{equation}
\label{eq:tv-H-definition}
    H=\sum_{i=1}^K\left(
       \sum_{m=2}^k a_{m,i}r_i^{-(m-1)/2}
          M_m\mathbin{\lrcorner}U_{m,i}
       +a_{k,i}R_{k,r_i}(W_i)\mathbin{\lrcorner}U_{k,i}\right),
\end{equation}
where $\lrcorner$ contracts the derivative indices and leaves the output index free.
Each integration by parts uses a fresh block, so no derivative falls on a score or on the remainder coefficient.

Independence and the bound we assume on second moments of scores give, for any tensor $T$ depending only on $W_i$,
\[
    \mathbb{E}[\|T\mathbin{\lrcorner}U_{m,i}\|_2^2]
    \le L^m\mathbb{E}[\|T\|_{\mathrm F}^2].
\]
Indeed, the second-moment matrix of $U_{m,i}$ is the tensor product of those of its score factors and is bounded by $L^mI_{N^m}$, independently of $T$.
This applies both to the deterministic tensors $M_m$ and to $R_{k,r_i}(W_i)$.
Since $r_i\asymp_k r$ and the coefficients $a_{m,i}$ are bounded in terms of $k$, the triangle inequality and \eqref{eq:tv-block-remainder} prove \eqref{eq:tv-H-bound}.
\end{proof}

\subsubsection{Score vector estimate}
\label{sec:tv-gaussian-score}

We now verify the score assumption of Lemma~\ref{lem:tv-block-representation} for Gaussian underlying vectors.
We first bound inverse moments of the empirical matrix that appears in the score calculation.

\begin{proposition}[Inverse moments of an empirical matrix]
\label{prop:tv-inverse-moments}
Let $N=\binom np$ and let $v^{(1)},\ldots,v^{(r)}$ be independent with law $\nu^{\otimes n}$, where $\nu$ is log-concave and $\Ex_{v\sim\nu}[\exp(cv^2)]<\infty$ for some $c>0$.
Let $B(v)$ be an $N\times N$ matrix whose entries are polynomials of degree at most $2p$, with
\begin{equation}
\label{eq:tv-single-jacobian}
    \mathbb{E}[B(v)]=I_N,\qquad
    0\preceq B(v)\preceq C\|v\|^{2p-2}I_N
\end{equation}
on the support of $\nu^{\otimes n}$, uniformly in $n$.
For $\Gamma=r^{-1}\sum_{j=1}^r B(v^{(j)})$ and each fixed $s>0$, for some constant $C_s > 0$, if $r\ge C_sn^p$ then we have
\begin{equation}
\label{eq:tv-inverse-moments}
    \mathbb{E}[\|\Gamma^{-1}\|_{\mathrm{op}}^s]\le C_s.
\end{equation}
\end{proposition}
\begin{proof}
For every unit vector $u$, $u^\top B(v)u$ is a nonnegative polynomial of degree at most $2p$ and mean 1.
The Carbery--Wright inequality \cite[Theorem~8]{CW-2001-ConvexBodyAnticoncentration} gives, for $t \in (0, 1)$,
\[
    \Prx[u^\top B(v)u\le t]\lesssim_p t^{1/(2p)}.
\]
Choose a $\varepsilon_t$-net $\mathcal U_t$ of the unit sphere in $\mathbb R^N$ with $\varepsilon_t=t^2/(4n^{p-1})$.
We may choose $|\mathcal U_t|\le(Cn^{p-1}/t^2)^N$.
If $\|\Gamma\|_{\mathrm{op}}\le n^{p-1}/t$ and $\lambda_{\min}(\Gamma)\le t$, approximating a unit eigenvector for $\lambda_{\min}(\Gamma)$ by a vector $u\in\mathcal U_t$ gives
\[
    u^\top\Gamma u\le t+2\|\Gamma\|_{\mathrm{op}}\varepsilon_t\le2t.
\]
At least half the summands of that quadratic form must then be at most $4t$.
Independence and a union bound over the net and these subsets of summands give, for any fixed $q>0$,
\begin{equation}
\label{eq:tv-small-eigenvalue}
    \Prx[\lambda_{\min}(\Gamma)\le t]
    \le C_qt^q+(Cn^{p-1})^N C^r t^{r/(4p)-2N}.
\end{equation}
For the first term, \eqref{eq:tv-single-jacobian} and the moment assumption imply
$\mathbb{E}[\|\Gamma\|_{\mathrm{op}}^q]\lesssim_q n^{q(p-1)}$.

Choose a sufficiently large constant $L$ and truncate each $B(v)$ to $\|v\|^2\le Ln$.
The truncated matrix has norm $\lesssim n^{p-1}$ and expectation at least $I_N/2$.
Indeed, the discarded expectation has norm at most a polynomial in $n$ times $\exp(-cn)$, with the exponential rate made arbitrarily large by increasing $L$.
The matrix Chernoff inequality \cite[Theorem~5.1.1]{Tropp-2015-MatrixConc} gives
\begin{equation}
\label{eq:tv-fixed-eigenvalue}
    \Prx[\lambda_{\min}(\Gamma)\le1/4]
    \le N\exp(-cr/n^{p-1}).
\end{equation}
Integrating \eqref{eq:tv-small-eigenvalue} against $s t^{-s-1}$ for $0<t<\delta/n$ and taking $q>s$, the second term contributes at most
\[
    \frac{C^{r+N}\delta^{r/(4p)-2N-s}
    n^{(p+1)N-r/(4p)+s}}{r/(4p)-2N-s}.
\]
For $\delta/n\le t\le1/4$, \eqref{eq:tv-fixed-eigenvalue} contributes at most
$C_s n^sN\exp(-cr/n^{p-1})$.
These quantities are uniformly bounded when $r\ge C_sn^p$ and then $\delta>0$ is sufficiently small.
This proves \eqref{eq:tv-inverse-moments}.
\end{proof}

\begin{lemma}[Score bound for Gaussian underlying vectors]
\label{lem:tv-block-score}
For $\nu=\mathcal N(0,1)$ and $r\ge Cn^p$, $Y_r$ has a density and a square-integrable weak score $\rho$, in the sense of Lemma~\ref{lem:stein-error-distance}, satisfying
\begin{equation}
\label{eq:tv-score-bound}
    \mathbb{E}[\rho(Y_r)\rho(Y_r)^\top]\preceq CI_N.
\end{equation}
\end{lemma}
\begin{proof}
Write $G=(v^{(1)},\ldots,v^{(r)})\in\mathbb R^{nr}$ for the underlying Gaussian vector, and let $J=D_GY_r$ be the $N\times nr$ Jacobian of the observation map.
We write $\partial_t$ for differentiation with respect to $G_t$.
Put
\[
    \Sigma=JJ^\top,\qquad
    \frac{\Sigma}p=\frac1r\sum_{j=1}^r B(v^{(j)}),\qquad
    B(v)=\frac1p DX(v)DX(v)^\top.
\]
Here $DX$ is the ordinary Jacobian of the map $X$.
Centering and independence give $\mathbb{E}[B(v)]=I_N$.
For a unit vector $u$, Cauchy--Schwarz in each sum
$\sum_{A\ni i}u_A\prod_{j\in A\setminus\{i\}}v_j$, followed by summation over $i$, gives the norm bound in \eqref{eq:tv-single-jacobian}.
Proposition~\ref{prop:tv-inverse-moments} therefore bounds each required fixed inverse moment of $\Sigma$ for $r\ge Cn^p$, with $C$ chosen for those moment orders.
In particular, $J$ has full row rank almost surely, so the coarea formula \cite[Section~3.4]{EG-2015-MeasureTheoryFineProperties} gives a density for $Y_r$.

Fix a unit output direction $u\in\mathbb R^N$ and set
\[
    w=J^\top\Sigma^{-1}u,\qquad U=\|\Sigma^{-1}\|_{\mathrm{op}},\qquad
    P=J^\top\Sigma^{-1}J.
\]
The choice of $w$ ensures $Jw=u$.
Thus, for $h\in C_c^\infty(\mathbb R^N)$, the chain rule gives
\[
    \langle\nabla_G(h(Y_r)),w\rangle
    =\langle\nabla h(Y_r),Jw\rangle
    =\langle u,\nabla h(Y_r)\rangle.
\]
This expresses a derivative in the observation direction $u$ as a derivative in the underlying Gaussian variables, where we can integrate by parts.
We first bound $w$ and $Dw$ to control the second moment of the coefficient produced by that integration.
The matrix $P$ is an orthogonal projection.
Differentiating $w$ gives
\[
    \partial_t w=(I-P)(\partial_tJ)^\top\Sigma^{-1}u
       -J^\top\Sigma^{-1}(\partial_tJ)w.
\]
Consequently, for
\[
    S=\sum_t(\partial_tJ)(\partial_tJ)^\top,\qquad
    Q=\sum_t(\partial_tJ)^\top(\partial_tJ),
\]
we have
\begin{equation}
\label{eq:tv-lift-energy}
    \|w\|^2\le U,\qquad
    \|Dw\|_{\mathrm{F}}^2
    \le2U^2\left(\|S\|_{\mathrm{op}}+\|Q\|_{\mathrm{op}}\right).
\end{equation}
Here $Dw$ is the Jacobian of $w: \mathbb{R}^{nr} \to \mathbb{R}^{nr}$.

We next bound the second moments of $\|S\|_{\mathrm{op}}$ and $\|Q\|_{\mathrm{op}}$.
For $i\in[n]$, let $D_i(v)=\partial_i DX(v)$, an $N\times n$ matrix with entries
\[
    (D_i(v))_{A,j}
    =\bm 1\{i\ne j,\ \{i,j\}\subseteq A\}
      \prod_{\ell\in A\setminus\{i,j\}}v_\ell.
\]
Put $S_0(v)=\sum_iD_i(v)D_i(v)^\top$ and $Q_0(v)=\sum_iD_i(v)^\top D_i(v)$.
Since the summands use disjoint input coordinates,
\[
    S=\frac1r\sum_{a=1}^r S_0(v^{(a)}),\qquad
    Q=\frac1r\operatorname{diag}\bigl(Q_0(v^{(1)}),\ldots,Q_0(v^{(r)})\bigr).
\]
For $u\in\mathbb R^N$, Cauchy--Schwarz gives
\begin{align*}
    u^\top S_0(v)u
    &=\sum_{i\ne j}\left(
       \sum_{A\supseteq\{i,j\}}u_A
       \prod_{\ell\in A\setminus\{i,j\}}v_\ell\right)^2\\
    &\le\frac{\|v\|^{2p-4}}{(p-2)!}
       \sum_{i\ne j}\sum_{A\supseteq\{i,j\}}u_A^2
     =\frac{p(p-1)}{(p-2)!}\|v\|^{2p-4}\|u\|^2.
\end{align*}
For $A\ne B$, every term in $(S_0(v))_{A,B}$ contains a centered coordinate to the first power, so its expectation is zero.
Each diagonal entry has $p(p-1)$ terms with expectation 1.
Thus
\begin{equation}
\label{eq:tv-hessian-estimates}
    \mathbb{E}[S_0(v)]=p(p-1)I_N,\qquad
    \|S_0(v)\|_{\mathrm{op}}\lesssim_p\|v\|^{2p-4}.
\end{equation}

Fix a sufficiently large constant $L$ and define
\[
    S_{\mathrm{trunc}} \colonequals \frac1r\sum_{a=1}^r
      S_0(v^{(a)})\bm 1\{\|v^{(a)}\|^2\le Ln\}.
\]
The matrices before division by $r$ are independent, positive semidefinite, and bounded above by $Cn^{p-2}I_N$, while $\mathbb{E}[S_{\mathrm{trunc}}]\preceq p(p-1)I_N$.
The upper-tail matrix Chernoff bound \cite[Theorem~5.1.1, Eq.~(5.1.6)]{Tropp-2015-MatrixConc} therefore gives, for a sufficiently large constant $t_0$ and $t \geq t_0$,
\[
    \Prx[\|S_{\mathrm{trunc}}\|_{\mathrm{op}}>t]
    \le N\exp(-crt/n^{p-2}).
\]
Integrating this tail with $b=cr/n^{p-2}$ gives
\begin{align*}
    \mathbb{E}[\|S_{\mathrm{trunc}}\|_{\mathrm{op}}^2]
    &\le t_0^2+2N\int_{t_0}^\infty t\exp(-bt)\,dt\\
    &=t_0^2+2N\exp(-bt_0)(t_0/b+1/b^2)\\
    &\lesssim_p 1,
\end{align*}
where the last inequality uses $r\ge Cn^p$ and $N\le n^p$.
Minkowski's inequality and the Gaussian tail bound give
\[
    \bigl\|\|S-S_{\mathrm{trunc}}\|_{\mathrm{op}}\bigr\|_{L^2}
    \lesssim_p\left(\mathbb{E}[\|v\|^{4p-8}\bm 1\{\|v\|^2>Ln\}]\right)^{1/2}
    \lesssim_p n^{p-2}\exp(-cn).
\]
Together these bounds control the second moment of $\|S\|_{\mathrm{op}}$ uniformly for all $r\ge Cn^p$.

For $z\in\mathbb R^n$, another application of Cauchy--Schwarz gives
\begin{align*}
    z^\top Q_0(v)z
    &=\sum_i\sum_{A\ni i}\left(
      \sum_{j\in A\setminus\{i\}}z_j
      \prod_{\ell\in A\setminus\{i,j\}}v_\ell\right)^2\\
    &\le(p-1)\sum_j z_j^2\sum_{i\ne j}
      \sum_{A\supseteq\{i,j\}}
      \prod_{\ell\in A\setminus\{i,j\}}v_\ell^2 \\
  &\lesssim_p n\|v\|^{2p-4}\|z\|^2.
\end{align*}
Thus $\mathbb{E}[\|Q_0(v)\|_{\mathrm{op}}^{2p}]\lesssim_p n^{2p(p-1)}$ by Gaussian moment bounds.
Using the block diagonal form of $Q$ and bounding a maximum by a sum of $2p$-th powers, we obtain
\[
    \bigl\|\|Q\|_{\mathrm{op}}\bigr\|_{L^2}
    \le\frac1r\left(r\mathbb{E}[\|Q_0(v)\|_{\mathrm{op}}^{2p}]\right)^{1/(2p)}
    \lesssim_p n^{p-1}r^{-1+1/(2p)}.
\]
In particular, for $r\ge Cn^p$,
\begin{equation}
\label{eq:tv-hessian-moments}
    \bigl\|\|S\|_{\mathrm{op}}\bigr\|_{L^2}\lesssim_p 1,\qquad
    \bigl\|\|Q\|_{\mathrm{op}}\bigr\|_{L^2}\lesssim_p 1.
\end{equation}

We now use these bounds to construct and bound the weak score.
Gaussian integration by parts in $G$ gives
\[
    \mathbb{E}[\langle\nabla_G g(G),w(G)\rangle]
    =\mathbb{E}[g(G)\delta w(G)],
    \qquad
    \delta w=\sum_{t=1}^{nr}(G_tw_t-\partial_t w_t).
\]
This identity says that the operator $\delta$ is the adjoint of the gradient with respect to Gaussian expectation \cite[Section~1.3]{Nualart-2006-MalliavinCalculus}.
We apply it with $g(G)=h(Y_r)$, using the chain-rule identity above to obtain the derivative in direction $u$.
To control the resulting coefficient $\delta w$, two further integrations by parts give
\[
    \mathbb{E}[(\delta w)^2]
    =\mathbb{E}\left[\|w\|^2+\sum_{t,s=1}^{nr}(\partial_t w_s)(\partial_s w_t)\right]
    \le\mathbb{E}[\|w\|^2+\|Dw\|_{\mathrm{F}}^2].
\]
Writing $U=\|\Sigma^{-1}\|_{\mathrm{op}}$ as above, \eqref{eq:tv-lift-energy} and Cauchy--Schwarz bound the last expression by
\[
    \mathbb{E}[U]
    +2\mathbb{E}[U^4]^{1/2}
      \left(\bigl\|\|S\|_{\mathrm{op}}\bigr\|_{L^2}
           +\bigl\|\|Q\|_{\mathrm{op}}\bigr\|_{L^2}\right)
    \lesssim_p 1,
\]
using \eqref{eq:tv-inverse-moments} with moment orders 1 and 4, and \eqref{eq:tv-hessian-moments}.
Since $w=J^\top\Sigma^{-1}u$ depends linearly on $u$, the identities
$\langle u,\rho(Y_r)\rangle=-\mathbb{E}[\delta w\mid Y_r]$ define a vector $\rho(Y_r)$.
The integration-by-parts identity just established then gives, for $h\in C_c^\infty(\mathbb R^N)$,
\[
    \mathbb{E}[\langle u,\nabla h(Y_r)\rangle]
    =\mathbb{E}[h(Y_r)\delta w]
    =-\mathbb{E}[h(Y_r)\langle u,\rho(Y_r)\rangle].
\]
Thus $\rho$ is a weak score, and a conditional Jensen inequality gives, for every unit $u$,
\[
    u^\top\mathbb{E}[\rho(Y_r)\rho(Y_r)^\top]u
    =\mathbb{E}\left[\mathbb{E}[\delta w\mid Y_r]^2\right]
    \le\mathbb{E}[(\delta w)^2]\lesssim_p 1,
\]
which proves \eqref{eq:tv-score-bound}.
To justify these integrations, first use $w_\eta=J^\top(\Sigma+\eta I_N)^{-1}u$ with smooth compactly supported cutoffs.
The bounds in \eqref{eq:tv-lift-energy} hold for $w_\eta,Dw_\eta$ uniformly in $\eta>0$, and their right-hand sides are integrable by the estimates above.
After removing the cutoffs, dominated convergence gives $w_\eta\to w$ and $Dw_\eta\to Dw$ in $L^2$ as $\eta\downarrow0$.
The displayed estimate for $\delta$ also gives $\delta w_\eta\to\delta w$ in $L^2$, so the identities pass to the limit.
\end{proof}

\subsubsection{Bounds on cumulants and Stein identity error}
\label{sec:tv-cumulants}

It remains to bound the coefficient tensors in Lemma~\ref{lem:tv-block-representation}.
The required cumulant estimates follow by counting overlaps among tensor coordinates.

\begin{lemma}[Cumulant bounds for one summand]
\label{lem:tv-cumulants}
For any $\nu$ of mean zero and variance one with moments of all orders and every fixed $m\ge2$,
\begin{equation}
\label{eq:tv-cumulant-bound}
    \|\kappa_m(X)\|_{\mathrm F}
    \lesssim_{p,m,\nu}n^{\lfloor mp/2\rfloor/2}.
\end{equation}
For $\nu = \mathcal{N}(0, 1)$ in particular, we have
\begin{equation}
\label{eq:tv-gaussian-third-cumulant}
    \|\kappa_3(X)\|_{\mathrm F}^2
    =\begin{cases}
        \dfrac{(n)_{3p/2}}{((p/2)!)^3},&p\text{ even},\\
        0,&p\text{ odd},
    \end{cases}
\end{equation}
where $(n)_j=n(n-1)\cdots(n-j+1)$, with $(n)_j=0$ for $j>n$.
\end{lemma}
\begin{proof}
If the hypergraph with edges $A_1,\ldots,A_m$ is disconnected, the variables $X_{A_1},\ldots,X_{A_m}$ split into independent groups, so Proposition~\ref{prop:cumulant-independence} makes their joint cumulant zero.
If a vector coordinate occurs in exactly one of these sets, every term in the cumulant formula of Proposition~\ref{prop:moment-cumulant-identities} vanishes by centering and independence.
Thus a nonzero entry of $\kappa_m(X)$ uses at most $\lfloor mp/2\rfloor$ distinct vector coordinates.
There are only finitely many overlap patterns for fixed $p,m$, and each entry is bounded in terms of $p,m,\nu$.
Counting the realizations of these patterns proves \eqref{eq:tv-cumulant-bound}.

For Gaussian inputs, $\kappa_3(X)_{A,B,C}=\mathbb{E}[X_AX_BX_C]$ is nonzero exactly when every vector coordinate occurs twice.
Then the three pairwise intersections are disjoint and each has size $p/2$, and the moment equals one.
Counting the ordered triples of these intersections gives \eqref{eq:tv-gaussian-third-cumulant}.
\end{proof}

\begin{proposition}[Quantitative Stein identity error bound]
\label{prop:tv-block-reduction}
For Gaussian inputs and $r\ge C_p n^p$, the representation~\eqref{eq:tv-H-representation} holds with
\begin{equation}
\label{eq:tv-final-fisher}
    I(\mathrm{Law}(Y_r)\mid\mathcal N)^{1/2}
    \le\|H\|_{L^2}
    \lesssim_p\left(
        \frac{\|\kappa_3(X)\|_{\mathrm F}}{\sqrt r}
        +\frac{n^p}{r}\right).
\end{equation}
\end{proposition}
\begin{proof}
Take $k=p$ and $K=p+1$ in Lemma~\ref{lem:tv-block-representation}.
For a sufficiently large constant in $r\ge C_p n^p$, every block size $r_i$ satisfies the hypothesis of Lemma~\ref{lem:tv-block-score}.
That lemma supplies the score bound, while Lemma~\ref{lem:tv-cumulant-expansion} supplies the expansion and remainder estimate.
Therefore \eqref{eq:tv-H-bound} gives
\[
    \|H\|_{L^2}
    \lesssim_p\left(
       \sum_{m=2}^p\frac{\|M_m\|_{\mathrm F}}{r^{(m-1)/2}}
       +\frac{n^{p^2/2}}{r^{p/2}}\right).
\]
By Proposition~\ref{prop:tv-kernels} and Lemma~\ref{lem:tv-cumulants}, for $3\le m\le p$ and $r\ge n^p$,
\[
    \frac{\|M_m\|_{\mathrm F}}{r^{(m-1)/2}}
    \lesssim_p\frac{n^{\lfloor(m+1)p/2\rfloor/2}}{r^{(m-1)/2}}
    \le\frac{n^p}{r}.
\]
The remainder obeys
\[
    \frac{n^{p^2/2}}{r^{p/2}}
    =\left(\frac{n^p}{r}\right)^{p/2}
    \le\frac{n^p}{r}.
\]
Together with $M_2=\kappa_3(X)/2$, these estimates prove the bound on $H$.
Lemma~\ref{lem:tv-block-score} also supplies a weak score for $Y_r$ itself.
The representation~\eqref{eq:tv-H-representation} and Lemma~\ref{lem:stein-error-distance} then give $I(\Law(Y_r)\mid\mathcal N)^{1/2}\le\|H\|_{L^2}$.
\end{proof}

\subsubsection{Full tensors: Proof of Theorem~\ref{thm:full-tensor-tv}}
\label{sec:pf:thm:full-tensor-tv}

\begin{proof}[Proof of Theorem~\ref{thm:full-tensor-tv}]
For $r\ge C_p n^p$, Proposition~\ref{prop:tv-block-reduction} and Lemma~\ref{lem:stein-error-distance} give
\begin{equation}
\label{eq:tv-block-conclusion}
    d_{\mathrm{TV}}(\Law(Y_r),\mathcal N)
    \lesssim_p\min\left\{1,
        \frac{\|\kappa_3(X)\|_{\mathrm F}}{\sqrt r}
        +\frac{n^p}{r}\right\}.
\end{equation}
By Lemma~\ref{lem:tv-cumulants}, $\|\kappa_3(X)\|_{\mathrm F}\lesssim_p n^{3p/4}$ for even $p$ and $\kappa_3(X)=0$ for odd $p$.
For even $p$, absorb $n^p/r$ into $n^{3p/4}/\sqrt r$ to obtain the stated rate.
For odd $p$, only $n^p/r$ remains.
The range $r<C_p n^p$ follows from $d_{\mathrm{TV}}\le1$ after increasing the final constant.
The Wasserstein bound mentioned before the theorem follows from the same two results for $r\ge C_p n^p$.
\end{proof}

\subsubsection{Diluted models: Proof of Theorem~\ref{thm:diluted-tensor-tv}}
\label{sec:pf:thm:diluted-tensor-tv}

For $c\ge1/2$, the undiluted bound and contraction of TV distance under the Ornstein--Uhlenbeck semigroup suffice.
For $c<1/2$, we reuse the cumulant expansion and transfer derivatives directly to the added Gaussian noise, avoiding the score function estimates used above.

\begin{proof}[Proof of Theorem~\ref{thm:diluted-tensor-tv}]
Write $a=c^2n^{p-1}$ and $b=r_{\comp}(n;c)$, so $r_{\stat}(n;c)=\max\{a,b\}$.
For $c\ge1/2$, $b\asymp_p r_*(n)$, so Theorem~\ref{thm:full-tensor-tv} proves the claim: applying the same Gaussian dilution to both laws contracts TV distance and preserves the null model.
It remains to consider $0<c<1/2$.
Set $\sigma=\sqrt{1-c^2}$ and $Z=cY_r+\sigma G$, where $G\sim\mathcal N$ is independent of $Y_r$.
Gaussian smoothing gives $Z$ a density and a square-integrable weak score $\rho_Z(Z)=-\sigma^{-1}\mathbb E[G\mid Z]$.
Thus Lemma~\ref{lem:stein-error-distance} applies to $Z$.

We first integrate by parts in $G$, and then apply Lemma~\ref{lem:tv-cumulant-expansion} conditional on $G$.
The chain rule gives, for every fixed $k\ge3$,
\begin{align*}
    \mathcal E_Z(f)
    &=\sum_{m=2}^k c^{m+1}r^{-(m-1)/2}
        \langle M_m,\mathbb E[D^mf(Z)]\rangle\\
    &\hspace{1cm}+c^{k+1}\mathbb E[\langle R_{k,r}(Y_r),D^kf(Z)\rangle].
\end{align*}
We transfer these derivatives directly to $G$ by Gaussian integration by parts.
Let $\mathcal H_m(G)$ be the order-$m$ Hermite tensor, whose entry at $(i_1,\ldots,i_m)$ is $H_\alpha(G)$ for the multi-index $\alpha$ recording these indices' multiplicities.
Transferring $m$ derivatives introduces $\sigma^{-m}\mathcal H_m(G)$.
For any tensor $T$ independent of $G$, orthogonality of Hermite polynomials gives
\[
    \|T\mathbin{\lrcorner}\mathcal H_m(G)\|_{L^2}
    \le\sqrt{m!}\,\|T\|_{L^2(\mathrm F)}.
\]
This applies to both $M_m$ and $R_{k,r}(Y_r)$.
Since $\sigma\ge\sqrt3/2$, Cauchy--Schwarz and Lemma~\ref{lem:stein-error-distance} now give
\begin{equation}
\label{eq:tv-diluted-expansion-bound}
    d_{\mathrm{TV}}(\Law(Z),\mathcal N)
    \lesssim_{p,k}\left(
       \sum_{m=2}^k
       \frac{c^{m+1}\|\kappa_{m+1}(X)\|_{\mathrm F}}{r^{(m-1)/2}}
       +c n^{p/2}\left(\frac ar\right)^{k/2}\right).
\end{equation}
Here we used the remainder bound in Lemma~\ref{lem:tv-cumulant-expansion}.

By Lemma~\ref{lem:tv-cumulants}, the $m=2$ term is $\lesssim_p\sqrt{b/r}$ for even $p$ and vanishes for odd $p$.
For $m\ge3$, we have
\[
    \frac{c^{m+1}\|\kappa_{m+1}(X)\|_{\mathrm F}}{r^{(m-1)/2}}
    \lesssim_{p,k}\frac{c^4n^p}{r}
      \left(\frac{c^2n^{p/2}}r\right)^{(m-3)/2}.
\]
Since $c^2n^{p/2}\le a$, $c^4n^p\le b$, and $r\ge n^{\delta}\max\{a,b\}$, \eqref{eq:tv-diluted-expansion-bound} gives
\[
    d_{\mathrm{TV}}(\Law(Z),\mathcal N)
    \lesssim_{p,k}\left(n^{-\delta/2}+n^{-\delta}
        +n^{(p-k\delta)/2}\right).
\]
Choosing a fixed $k\ge3$ with $k\delta>p$ then gives the result.
\end{proof}

\subsection{Computational upper bounds from strongly separating low-degree polynomials}
\label{sec:polynomial-tests}

We use the following parameters in the proofs of the polynomial tests and low-degree bounds.
\begin{definition}
\label{defin:q,s}
Define
\[
(q,s)
\colonequals
\begin{cases}
\displaystyle
\left(3,\frac{3p}{2}\right),
& p \text{ even},\\[3mm]
\displaystyle
\left(3,\frac{3p-1}{2}\right),
& p \text{ odd and } \tau\neq 0,\\[3mm]
(4,2p),
& p \text{ odd and } \tau=0.
\end{cases}
\]
Here $q$ is the first degree at which the cumulants of $\mu_n$ differ from those of $\mathcal N$, and $s$ is the maximal number of vertices in a corresponding shape.
We have $r_*(n)=n^{s/(q-2)}$.
When $\nu$ has $(K,C)$-bounded moments (Definition~\ref{def:bounded-moments}), fix such a $K=K(\nu)$ and take
\[
\gamma_{p,\nu}
\coloneqq
\frac{pq}{q-2}(2K+1).
\]
\end{definition}

\subsubsection{Full tensors: Proof of Theorem~\ref{thm:full-tensor-polynomial-tests}}
\label{sec:pf:thm:full-tensor-polynomial-tests}

\begin{proof}
Let $q,s$ be as in Definition~\ref{defin:q,s}.
All constants below may depend on $p,\nu$.
We note that while the cases $q = 3$ and $q = 4$ involve different test statistics, we treat them both at once below.

Choose a fixed $p$-uniform hypergraph $H$ with $q$ distinct edges and $s$ vertices.
For $q=3$, take pairwise disjoint sets $A,B,C,D$ with $|A|=|B|=|C|=\lfloor p/2\rfloor$ and $|D|=p-2\lfloor p/2\rfloor$, and set
\[
    E(H)=\{A\cup B\cup D,\ B\cup C\cup D,\ C\cup A\cup D\}.
\]
For $q=4$, take disjoint sets $A,B$ of size $p-1$ and distinct vertices $i,j$ outside $A\cup B$, and set
\[
    E(H)=\{A\cup\{i\},\ A\cup\{j\},\ B\cup\{i\},\ B\cup\{j\}\}.
\]
Writing $\mathcal I_H$ for the injective maps $V(H)\hookrightarrow[n]$, define
\[
    Z_\phi(y)\colonequals\prod_{e\in E(H)}y_{\phi(e)},
    \qquad T_H(y)\colonequals\sum_{\phi\in\mathcal I_H}Z_\phi(y).
\]
Since the edges are distinct, centering and isotropy leave only the partition with one block in the moment--cumulant expansion of $\Ex_{\mathbb P_n}[Z_\phi]$.
Propositions~\ref{prop:moment-cumulant-identities} and~\ref{prop:cumulant-independence} therefore give
\begin{align*}
    \Ex_{\mathbb P_n}[Z_\phi]
    &=r^{1-q/2}\kappa(X_{\phi(e)}:e\in E(H))
      =r^{1-q/2}\prod_{v\in V(H)}m_{\deg_H(v)}(\nu)
      =\theta_Hr^{1-q/2},\\
    \theta_H&\colonequals
    \begin{cases}
        \tau,&q=3\text{ and }p\text{ odd},\\
        1,&\text{otherwise}.
    \end{cases}
\end{align*}
Here we used cumulant additivity and homogeneity for the normalized sum, as in Proposition~\ref{prop:scaling-cumulants}.
Thus, with $(n)_s=n(n-1)\cdots(n-s+1)$,
\[
    \Ex_{\mathbb Q_n}[T_H]=0,
    \qquad M_H\colonequals\Ex_{\mathbb P_n}[T_H]
    =\theta_H(n)_s r^{1-q/2}.
\]

Under the Gaussian null model, independence gives
\[
    \Ex_{\mathbb Q_n}[Z_\phi Z_\psi]
    =\bm 1\{\phi(E(H))=\psi(E(H))\},
    \qquad \Varx_{\mathbb Q_n}[T_H]\le s!(n)_s.
\]
For the alternative, enumerate $E(H)=\{e_1,\ldots,e_q\}$ and, for each $\phi,\psi\in\mathcal I_H$, put
$F_j=\phi(e_j)$ and $F_{q+j}=\psi(e_j)$ for $j\in[q]$.
Let $\Part_\times([2q])$ consist of the partitions with at least one block meeting both $[q]$ and $\{q+1,\ldots,2q\}$.
Subtracting the two moment--cumulant expansions gives
\[
    \Covx_{\mathbb P_n}[Z_\phi,Z_\psi]
    =\sum_{\pi\in\Part_\times([2q])}
      r^{|\pi|-q}\prod_{B\in\pi}\kappa(X_{F_j}:j\in B).
\]
This is the covariance case of the Leonov--Shiryaev product formula \cite[Proposition~4.3]{Speed-1983-CumulantsPartitonLattices}.
Across two copies of $H$, each underlying coordinate occurs at most six times, so the finite sixth moment assumption bounds all cumulant products by a constant.

Singleton blocks vanish, and isotropy gives
\[
    \kappa(X_F,X_G)=\bm 1\{F=G\}.
\]
When $q=4$, third cumulants also vanish.
Indeed, for any $p$-sets $F,G,J$, define $d_v=\bm 1\{v\in F\}+\bm 1\{v\in G\}+\bm 1\{v\in J\}$.
Then
\[
    \kappa(X_F,X_G,X_J)=\prod_{v\in F\cup G\cup J}m_{d_v}(\nu)=0,
\]
because $\sum_vd_v=3p$ is odd and $m_1(\nu)=\tau=0$.

Write $S_\phi=\phi(V(H))$.
For each fixed $1\le\ell\le s$,
\[
    \#\{(\phi,\psi)\in\mathcal I_H^2:|S_\phi\cap S_\psi|\ge\ell\}
    \le Cn^{2s-\ell}.
\]
A nonzero block meeting both copies requires $S_\phi\cap S_\psi\ne\varnothing$, by Proposition~\ref{prop:cumulant-independence}.
Every nonzero pair block matches an edge from one copy to an identical edge in the other.
If all blocks are pairs, the two edge sets coincide.
When $q=4$, the remaining partition type has block sizes $(2,2,4)$.
The two matched edges are distinct, so their union has at least $p+1$ vertices.
The possible nonzero contributions are therefore bounded as follows, for each fixed partition $\pi$:
\[
\begin{array}{c|c|c}
    \text{partition type}
    &\text{number of contributing }(\phi,\psi)
    &r^{|\pi|-q}\\ \hline
    |\pi|\le2 & O(n^{2s-1}) & \le r^{2-q}\\
    \text{all pairs} & O(n^s) & 1\\
    (2,2,4)\text{ when }q=4 & O(n^{4p-(p+1)}) & r^{-1}
\end{array}
\]
There are only finitely many partitions of $[2q]$.
Summing the absolute values of their contributions gives
\[
    \Varx_{\mathbb P_n}[T_H]\le C
    \begin{cases}
        n^{2s-1}/r+n^s,&q=3,\\
        n^{4p-1}/r^2+n^{3p-1}/r+n^{2p},&q=4.
    \end{cases}
\]

Since $\theta_H\ne0$ and $(n)_s\asymp n^s$, the bounds above imply
\[
    \frac{\Varx_{\mathbb P_n}[T_H]+\Varx_{\mathbb Q_n}[T_H]}{M_H^2}
    \le C
    \begin{cases}
        n^{-1}+r/n^s,&q=3,\\
        n^{-1}+r/n^{p+1}+r^2/n^{2p},&q=4.
    \end{cases}
\]
This tends to zero when $r=o(r_*(n))$, so the test statistic $f_n\colonequals\operatorname{sgn}(\theta_H)T_H$ achieves strong separation.
\end{proof}

\subsubsection{Diluted models: Proof of Theorem~\ref{thm:diluted-polynomial-tests}}
\label{sec:pf:thm:diluted-polynomial-tests}

\begin{proof}
We use the statistics $T_H$ from the preceding proof, writing $T_3$ for the cubic statistic when $q=3$ and $T_4$ for the quartic statistic.
The four-edge construction also applies when $q=3$.
Write $\widetilde{\mathbb P}_n=G_c(\mathbb P_n)$.
Gaussian dilution preserves the first two cumulants and multiplies each cumulant of order $j\ge3$ by $c^j$.
Thus the null means are zero, while the alternative means are
\[
    M_3\colonequals\Ex_{\widetilde{\mathbb P}_n}[T_3]
       =c^3\theta_H(n)_s/\sqrt r,
    \qquad
    M_4\colonequals\Ex_{\widetilde{\mathbb P}_n}[T_4]
       =c^4(n)_{2p}/r.
\]

We use the same covariance expansion and overlap counts as in the preceding proof.
For a degree-$\ell$ statistic, a partition $\pi$ of its $2\ell$ factors with $b_2$ pair blocks now carries the factor $c^{2\ell-2b_2}r^{|\pi|-\ell}$.
Each nonzero pair block matches an edge between the two copies, forcing at least $p$ common vertices.
Two such blocks match distinct edges and force at least $p+1$ common vertices.
Every other nonzero covariance term requires at least one common vertex.
These counts give
\begin{align*}
    \Varx_{\widetilde{\mathbb P}_n}[T_3]
    &\le C_{p,\nu}\left(
       \frac{c^6n^{2s-1}}r
       +\frac{c^4n^{2s-p}}r+n^s\right),\\
    \Varx_{\widetilde{\mathbb P}_n}[T_4]
    &\le C_{p,\nu}\left(
       \frac{c^8n^{4p-1}}{r^2}
       +\frac{c^6n^{3p}}r
       +\frac{c^4n^{3p-1}}r+n^{2p}\right).
\end{align*}
In the quartic bound, partitions of types $(2,6)$ and $(2,3,3)$ contribute to the second term, and type $(2,2,4)$ contributes to the third.
This includes the possible nonzero third cumulants when $q=3$.
The null variance bounds remain $O(n^s)$ and $O(n^{2p})$.
Dividing by the squared means gives
\begin{align*}
    \frac{\Varx_{\widetilde{\mathbb P}_n}[T_3]
       +\Varx_{\mathbb Q_n}[T_3]}{M_3^2}
    &\le C_{p,\nu}\left(
       \frac1n+\frac1{c^2n^p}+\frac r{c^6n^s}\right),\\
    \frac{\Varx_{\widetilde{\mathbb P}_n}[T_4]
       +\Varx_{\mathbb Q_n}[T_4]}{M_4^2}
    &\le C_{p,\nu}\left(
       \frac1n+\frac r{c^2n^p}
       +\frac r{c^4n^{p+1}}+\frac{r^2}{c^8n^{2p}}\right).
\end{align*}
If $r=o(c^6n^s)$, then $c^6n^s\to\infty$ and $c^2n^p=(c^6n^s)^{1/3}n^{p-s/3}\to\infty$, so $\operatorname{sgn}(\theta_H)T_3$ strongly separates.
If $r=o(c^4n^p)$, the quartic bound also tends to zero, since $c\le1$, so $T_4$ strongly separates.
The evaluation costs are the numbers of embeddings, $O(n^s)$ and $O(n^{2p})$.
Choosing the statistic corresponding to the larger of the two scales proves the claim for $r=o(r_{\comp}(n;c))$.
\end{proof}

\subsection{Computational lower bounds from low-degree polynomials}

\subsubsection{Full tensors: Proof of Theorem~\ref{thm:low-degree-1}}

We give a more specific statement involving the low-degree $\chi^2$-divergence.
\begin{theorem}[Low-degree $\chi^2$-divergence bounds for full tensors]
\label{thm:low-degree-1-chi-squared}
Let $(\mathbb Q_n,\mathbb P_n)$ be the Wigner--Wishart tensor testing
problem of Definition~\ref{def:wigner-wishart}, with parameters
$p$, $\nu$, and $r=r(n)$, where $p$ and $\nu$ are fixed as $n\to\infty$.
Assume that $\nu$ has $(K,C)$-bounded moments.
Let $q$, $s$, $r_*(n)$, and $\gamma_{p,\nu}$ be as in
Definition~\ref{defin:q,s}.
For a sufficiently small constant $\kappa=\kappa(p,K,C)>0$, if
\[
    r=\omega \left(r_*(n)\right), \qquad  D \le\kappa\left(
        \frac{r}{r_*(n)}
    \right)^{1/\gamma_{p,\nu}},
\]
then
\[
    \chi^2_{\leq D}
    \left(\mathbb P_n \,\middle|\, \mathbb Q_n\right)
    = O(1).
\]
Uniformly in $r\ge1$, we also have
\[
    \chi^2_{\le q}(\mathbb P_n\mid\mathbb Q_n)
    \asymp_{p,\nu}\frac{n^s}{r^{q-2}}
    =\left(\frac{r_*(n)}r\right)^{q-2}.
\]
Here $q\in\{3,4\}$ is the first degree at which the cumulants of $\mu_n$ differ from those of $\mathcal N$, as specified in Definition~\ref{defin:q,s}.
In particular, if $r=o(r_*(n))$, then
\[
    \chi^2_{\le 4}(\mathbb P_n\mid\mathbb Q_n)\longrightarrow\infty.
\]
\end{theorem}

\begin{proof}

Write $r_*=r_*(n)=n^{s/(q-2)}$, and let $\mathcal B=\mathcal A$ be the admissible class defined in Section~\ref{sec:masked-models}.
Its elements have degree at least $q$. In the notation of
Definition \ref{def:decomp-multigraph}, $f_{\mathcal B}(\alpha)=f(H_\alpha)$.

By Lemma \ref{lem:tensor-moment-bound}, $\mu_n$ has $(pK,C^p)$-bounded moments.
Thus Lemma \ref{lem:clt-cumulant-support-bound} gives
\[
    \chi^2_{\le D}(\mathbb P_n\mid\mathbb Q_n)
    \le
    \sum_{d=q}^{D}
    \left(\frac{16C^{2p}D^{2pK}}{r}\right)^d
    \sum_{\alpha\in\mathcal A_d}r^{2f(H_\alpha)},
\]
where $\mathcal A_d$ is the degree-$d$ decomposable class
defined in Definition \ref{def:decomp-class}. If $D<q$, the divergence
is zero. For any admissible shape $J$ with $e$ edges, the
minimum degree condition and Definition \ref{defin:q,s} imply
\[
    |V(J)|
    \le \left\lfloor\frac{pe}{2}\right\rfloor
    \le \frac{s}{q-2}(e-2).
\]
Applying this bound to a decomposition of a degree-$d$
shape $H$ into $f(H)$ admissible parts yields
\[
    |V(H)|\le\frac{s}{q-2}\bigl(d-2f(H)\bigr).
\]
Since $f(H)\le d/q$ and $r/r_*\to\infty$, for sufficiently large $n$ we have
\[
    n^{|V(H)|}r^{2f(H)-d}
    \le
    \left(\frac{r_*}{r}\right)^{d-2f(H)}
    \le
    \left(\frac{r_*}{r}\right)^{(q-2)d/q}.
\]

There are at most $(pd)^{pd}$ degree-$d$ shapes:
each is represented by a partition of the $pd$ vertex
positions in an ordered list of edges. Each shape $H$
has at most $n^{|V(H)|}$ realizations. Grouping by shape
and using $d\le D$, we therefore obtain
\begin{align*}
    \chi^2_{\le D}(\mathbb P_n\mid\mathbb Q_n)
    &\le
    \sum_{d=q}^{D}
    \left[
        16C^{2p}p^pD^{2pK+p}
        \left(\frac{r_*}{r}\right)^{(q-2)/q}
    \right]^d \\
    &\le
    \sum_{d=q}^{D}
    \left[
        16C^{2p}p^p
        \left( D\left(\frac{r_*}{r}\right)^{1/\gamma_{p,\nu}}\right)^{2pK+p}
    \right]^d
\end{align*}
We note that $\gamma_{p,\nu}=pq(2K+1)/(q-2)$.
Under the assumed degree restriction, choosing $\kappa$ sufficiently small makes the quantity in brackets at most $1/2$.
The resulting
geometric series is bounded, proving the first claim.

For the degree-$q$ estimate, the Gaussian-null identity of Theorem~\ref{thm:gaussian-null} gives
\[
    \chi^2_{\le q}(\mathbb P_n\mid\mathbb Q_n)
    =r^{2-q}\sum_{|\alpha|=q}\frac{\kappa_\alpha(\mu_n)^2}{\alpha!},
\]
since all cumulants of degrees below $q$ agree with those of the Gaussian.
The degree-three and degree-four counts at the end of the proof of Theorem~\ref{thm:gaussian-dilution-chi-squared} show that this sum is of order $n^s$.
This proves the estimate and its divergence when $r=o(r_*(n))$.
\end{proof}

\begin{proof}[Proof of Theorem~\ref{thm:low-degree-1}]
The claim follows from Theorem~\ref{thm:low-degree-1-chi-squared}
and Proposition~\ref{prop:chi-squared-testing}.
\end{proof}

\subsubsection{Diluted models: Proof of Theorem~\ref{thm:gaussian-dilution}}

\begin{theorem}[Low-degree $\chi^2$-divergence bounds for Gaussian dilution]
\label{thm:gaussian-dilution-chi-squared}
Let $(\mathbb Q_n,\mathbb P_n)$ be the Wigner--Wishart tensor testing
problem of Definition~\ref{def:wigner-wishart}, with $p$ and $\nu$ fixed as $n\to\infty$.
Assume that $\nu$ has $(K,C)$-bounded moments.
Let $q,s,\gamma_{p,\nu}$ be as in Definition~\ref{defin:q,s}.
For $c=c(n)\in(0,1]$, the scale in Theorem~\ref{thm:gaussian-dilution} can also be written as
\[
    r_{\comp}(n; c)
    =
    \max\left\{
        c^{2q/(q-2)}n^{s/(q-2)},
        c^4n^p
    \right\}.
\]
For a sufficiently small constant $\kappa=\kappa(p,K,C)>0$, if
\[
    r=\omega\bigl(r_{\comp}(n; c)\bigr),
    \qquad
    D\le\kappa
    \left(\frac{r}{r_{\comp}(n; c)}\right)^{1/\gamma_{p,\nu}},
\]
then
\[
    \chi^2_{\le D}
    \left(
        G_c(\mathbb P_n)
        \,\middle|\,
        \mathbb Q_n
    \right)
    =O(1).
\]
Uniformly in $r\ge1$ and $c\in(0,1]$, we also have
\[
    \chi^2_{\le4}(G_c(\mathbb P_n)\mid\mathbb Q_n)
    \asymp_{p,\nu}
    \begin{cases}
        \displaystyle \frac{c^6n^{3p/2}}r+\frac{c^8n^{2p}}{r^2},
        & p\text{ even},\\[2mm]
        \displaystyle \frac{c^6\tau^2n^{(3p-1)/2}}r+\frac{c^8n^{2p}}{r^2},
        & p\text{ odd}.
    \end{cases}
\]
In particular, if $r=o(r_{\comp}(n; c))$, then
\[
    \chi^2_{\le4}(G_c(\mathbb P_n)\mid\mathbb Q_n)\longrightarrow\infty.
\]
\end{theorem}

\begin{proof}
    Write $r_{\comp}=r_{\comp}(n; c)$. Gaussian dilution preserves the
first two cumulants and multiplies every cumulant of degree
$d\ge3$ by $c^d$. Thus, by Proposition \ref{prop:normalized-convolution-cumulants} and Lemma \ref{lem:tensor-moment-bound},
\[
    \left|\kappa_\alpha(G_c(\mathbb P_n))\right|
    \le
    r\left(\frac{2C^pc}{\sqrt r}\right)^{|\alpha|}
    (|\alpha|!)^{pK},
    \qquad \text{whenever } |\alpha|\ge3.
\]
Apply Lemma~\ref{lem:exponential-multiindex-bound} to $F=\psi_{G_c(\mathbb P_n)}-\psi_{\mathcal N}$ with $\delta=r$, $\varepsilon=2C^pc/\sqrt r$, and exponent $pK$ in place of $K$.
Together with Theorem~\ref{thm:gaussian-null}, this gives
\[
    \chi^2_{\le D}(G_c(\mathbb P_n)\mid\mathbb Q_n)
    \le
    \sum_{d=q}^D
    \left(\frac{16C^{2p}c^2D^{2pK}}r\right)^d
    \sum_{\alpha\in\mathcal A_d}r^{2f(H_\alpha)}.
\]

We claim that every admissible shape $J$ with $e$ edges satisfies,
for sufficiently large $n$,
\[
    c^{2e}n^{|V(J)|}r^{2-e}
    \le
    \left(\frac{r_{\comp}}{r}\right)^{(q-2)e/q}.
\]
If $e=3$, necessarily $q=3$, and this follows from
$|V(J)|\le s$ and $c^6n^s\le r_{\comp}$. If $e\ge4$, then
$|V(J)|\le pe/2$, so
\[
    c^{2e}n^{|V(J)|}r^{2-e}
    \le
    \left(\frac{c^4n^p}{r}\right)^{e/2}r^{2-e/2}
    \le
    \left(\frac{r_{\comp}}{r}\right)^{e/2}
    \le
    \left(\frac{r_{\comp}}{r}\right)^{(q-2)e/q},
\]
using $r\ge1$, $r_{\comp}/r<1$, and $q\in\{3,4\}$. Multiplying this bound over a decomposition of $H$ into
$f(H)$ admissible parts yields
\[
    c^{2d}n^{|V(H)|}r^{2f(H)-d}
    \le
    \left(\frac{r_{\comp}}{r}\right)^{(q-2)d/q}.
\]
The shape-counting argument from the proof of Theorem~\ref{thm:low-degree-1-chi-squared} therefore gives
\begin{align*}
    \chi^2_{\le D}(G_c(\mathbb P_n)\mid\mathbb Q_n)
    &\le
    \sum_{d=q}^D
    \left[
        16C^{2p}p^pD^{2pK+p}
        \left(\frac{r_{\comp}}{r}\right)^{(q-2)/q}
    \right]^d \\
    &\le
    \sum_{d=q}^{D}
    \left[
        16C^{2p}p^p
        \left( D\left(\frac{r_{\comp}}{r}\right)^{1/\gamma_{p,\nu}}\right)^{2pK+p}
    \right]^d.
\end{align*}

Since $\gamma_{p,\nu}=pq(2K+1)/(q-2)$, the assumed degree restriction makes the quantity in brackets at most $1/2$ when $\kappa$ is sufficiently small.
The resulting geometric series is bounded, proving the first claim.

For the degree-four estimate, write
\[
    A_j(n)\colonequals\sum_{|\alpha|=j}\frac{\kappa_\alpha(\mu_n)^2}{\alpha!},
    \qquad j\in\{3,4\}.
\]
The series $F=\psi_{G_c(\mathbb P_n)}-\psi_{\mathcal N}$ has only terms of degree at least three, so $\exp(F)=1+F$ through degree four.
Theorem~\ref{thm:gaussian-null} therefore gives the exact identity
\[
    \chi^2_{\le4}(G_c(\mathbb P_n)\mid\mathbb Q_n)
    =\frac{c^6}{r}A_3(n)+\frac{c^8}{r^2}A_4(n).
\]
Every contributing shape has minimum vertex degree two.
A three-edge shape thus has at most $\lfloor3p/2\rfloor$ vertices.
This maximum is attained with all vertex degrees two when $p$ is even, or with one degree-three vertex and all others of degree two when $p$ is odd.
The corresponding cumulants are $1$ and $\tau$, respectively; when $p$ is odd and $\tau=0$, all third cumulants vanish.
There are only finitely many shapes at each fixed degree, with cumulants bounded in terms of $p$ and $\nu$, so counting their embeddings gives the first contribution.

A four-edge shape has at most $2p$ vertices, giving $A_4(n)=O_{p,\nu}(n^{2p})$.
For a matching lower bound, take disjoint sets $S_1,S_2,S_3,S_4$ of sizes $1,p-1,1,p-1$ and edges $S_1\cup S_2,S_2\cup S_3,S_3\cup S_4,S_4\cup S_1$.
These give $\Theta_p(n^{2p})$ distinct configurations, each with fourth cumulant $1$: every vertex appears twice, and the four distinct tensor coordinates have zero pairwise covariance.
Thus $A_4(n)\asymp_{p,\nu}n^{2p}$, proving the estimate.
Its right-hand side diverges whenever $r=o(r_{\comp}(n; c))$, as claimed.
\end{proof}

\begin{proof}[Proof of Theorem~\ref{thm:gaussian-dilution}]
The claim follows from Theorem~\ref{thm:gaussian-dilution-chi-squared}
and Proposition~\ref{prop:chi-squared-testing}.
\end{proof}

\subsubsection{Masked models: Proof of Theorem~\ref{thm:masked-tensor-model}}

\begin{theorem}[Low-degree $\chi^2$-divergence bound for masked tensors]
\label{thm:masked-tensor-model-chi-squared}
Let $(\mathbb Q_n,\mathbb P_n)$ be the Wigner--Wishart tensor testing problem of Definition~\ref{def:wigner-wishart}.
Assume that $\nu$ has $(K,C)$-bounded moments.
For a sufficiently small constant $\kappa=\kappa(p,K,C)>0$, if
\[
    \zeta_{\leq D}(r; \Omega)=o(1), \qquad     D
    \le\kappa
    \left(
        \frac{1}{\zeta_{\leq D}(r; \Omega)}
    \right)^{1/(2pK+p)},
\]
then
\[
    \chi^2_{\leq D}\left(
        M_\Omega(\mathbb P_n)
        \,\middle|\,
        M_\Omega(\mathbb Q_n)
    \right)
    = O(1).
\]
\end{theorem}

\begin{proof}
Let $\eta\colonequals \zeta_{\leq D}(r; \Omega)$ and $\mathcal B_\Omega$ be the admissible multi-indices
supported on $\Omega$, and let $\mathcal A_d(\Omega)$
be the degree-$d$ multi-indices admitting a decomposition
into elements of $\mathcal B_\Omega$.
For such indices,
$f_{\mathcal B_\Omega}(\alpha)=f(H_\alpha)$.

Coordinate projection preserves centeredness, isotropy,
and the moment bound in Lemma~\ref{lem:tensor-moment-bound}. Moreover,
\[
    M_\Omega(\mathbb P_n)
    =
    r^{-1/2}\cdot M_\Omega(\mu_n)^{*r},
\]
and the vanishing criteria in Section~\ref{sec:tv-cumulants} imply that the cumulants of degree at least three of $M_\Omega(\mu_n)$ are supported on $\mathcal B_\Omega$. Thus Lemma~\ref{lem:clt-cumulant-support-bound} gives
\[
    \chi^2_{\le D}
    \bigl(M_\Omega(\mathbb P_n)\mid M_\Omega(\mathbb Q_n)\bigr)
    \le
    \sum_{d=q}^{D}
    \left(\frac{16C^{2p}D^{2pK}}r\right)^d
    \sum_{\alpha\in\mathcal A_d(\Omega)}
    r^{2f(H_\alpha)}.
\]

Fix a degree-$d$ shape $H$ admitting an admissible
decomposition, and choose one with the maximum number
$k=f(H)$ of parts:
\[
    E(H)=E(H_1)\sqcup\cdots\sqcup E(H_k).
\]
Every $H_j$ is irreducible, since otherwise refining that
part would contradict maximality. Put $d_j\colonequals |E(H_j)|$,
so that $q\le d_j\le D$ and $\sum_jd_j=d$.

Restricting an embedding of $H$ to the vertex sets of its
parts gives an injective map into the product of their
embedding sets. Consequently,
\[
    \rho(H,\Omega)
    \le
    \prod_{j=1}^{k}\rho(H_j,\Omega).
\]
By the definition of $\eta$, each irreducible part satisfies
\[
    \rho(H_j,\Omega)r^{2-d_j}\le\eta^{d_j}.
\]
Multiplying these inequalities yields
\[
    \rho(H,\Omega)r^{2f(H)-d}
    \le
    \prod_{j=1}^{k}
    \bigl(\rho(H_j,\Omega)r^{2-d_j}\bigr)
    \le\eta^d.
\]

There are at most $(pd)^{pd}$ degree-$d$ shapes, and the
number of realizations of a shape $H$ supported on $\Omega$
is at most $\rho(H,\Omega)$. Grouping by shape therefore
gives
\[
    \sum_{\alpha\in\mathcal A_d(\Omega)}
    r^{2f(H_\alpha)}
    \le
    (pd)^{pd}r^d\eta^d.
\]
Substituting into the preceding low-degree bound and
using $d\le D$, we obtain
\[
\begin{aligned}
    \chi^2_{\le D}
    \bigl(M_\Omega(\mathbb P_n)\mid M_\Omega(\mathbb Q_n)\bigr) &\le
    \sum_{d=q}^{D}
    \left[
        16C^{2p}p^pD^{2pK}d^{p}\eta
    \right]^d\\
    &\le
    \sum_{d=q}^{D}
    \left[
        16C^{2p}p^pD^{2pK+p}\eta
    \right]^d\\
    &\le \sum_{d=q}^{D}
    \left[
        16C^{2p}p^p (D\eta^{1/(2pK + p)})^{2pK +p}
    \right]^d.
\end{aligned}
\]
Under the assumed degree restriction
\[
    D\le\kappa\eta^{-1/(2pK+p)},
\]
choosing $\kappa$ sufficiently small makes the quantity in brackets at most $1/2$.
The resulting geometric series is bounded, proving
the claim.

If $\eta=0$, no irreducible admissible shape of degree
at most $D$ embeds into $\Omega$. The decomposition
argument above then implies
$\mathcal A_d(\Omega)=\varnothing$ for $q\le d\le D$,
so the low-degree divergence is zero.
\end{proof}

\begin{proof}[Proof of Theorem~\ref{thm:masked-tensor-model}]
The claim follows from Theorem~\ref{thm:masked-tensor-model-chi-squared}
and Proposition~\ref{prop:chi-squared-testing}.
\end{proof}

\subsection{Statistical upper bounds from brute-force search}
\label{sec:pf:thm:scan-detection}

\subsubsection{General analysis: Proof of Theorem~\ref{thm:scan-detection}}

\begin{proof}[Proof of Theorem~\ref{thm:scan-detection}]
Write $B_n\colonequals ca_n/\sqrt r$ and
\[
    \eta_n\colonequals\Prx_{X\sim\mu_n}\left[
        \max_{v\in\mathcal T_n}\langle v,X\rangle<a_n
    \right].
\]
The assumptions give $\eta_n\to0$ and
$B_n^2/(1+\log M_n)\to\infty$.
Under the null, each $\langle v,Y\rangle$ is a standard Gaussian.
A union bound therefore gives
\[
    \mathbb Q_n[T_n(Y)>B_n/2]
    \le M_n\exp(-B_n^2/8)\longrightarrow0.
\]

Under the alternative, write
\[
    Y=\frac c{\sqrt r}\sum_{i=1}^r X_i+\sqrt{1-c^2}\,G,
\]
where $X_1,\ldots,X_r\sim\mu_n$ and $G\sim\mathcal N$ are independent.
Choose $v_*(X_1)\in\mathcal T_n$ maximizing $\langle v,X_1\rangle$,
breaking ties by a fixed ordering of $\mathcal T_n$.
This choice is measurable and depends only on $X_1$.
Conditional on $X_1$,
\[
    \langle v_*(X_1),Y\rangle
    =\frac c{\sqrt r}\langle v_*(X_1),X_1\rangle+Z_n,
\]
where centeredness, isotropy, and independence give
\[
    \E[Z_n\mid X_1]=0,
    \qquad
    \Var[Z_n\mid X_1]
    =\frac{c^2(r-1)}r+1-c^2
    =1-\frac{c^2}{r}\le1.
\]
Except on an event of probability $\eta_n$, the first term is at least
$B_n$. Since $T_n(Y)\ge\langle v_*(X_1),Y\rangle$, a conditional
Chebyshev inequality gives
\[
    (G_c(\mathbb P_n))[T_n(Y)\le B_n/2]
    \le\eta_n+\frac4{B_n^2}\longrightarrow0.
\]
Together with the null bound, this proves strong detection.

For the discrete case, when $X\in\mathcal S_n$ the direction
$X/\|X\|_2$ belongs to $\mathcal T_n$ and has inner product
$\|X\|_2\ge b_n$ with $X$.
For the net case, when $X\in\mathcal S_n$, choose $v\in\mathcal T_n$ with
$\|v-X/\|X\|_2\|_2\le\varepsilon_n$.
Because both directions have unit norm,
\[
    \langle v,X\rangle
    =\|X\|_2\left(1-\frac12
        \left\|v-\frac X{\|X\|_2}\right\|_2^2\right)
    \ge(1-\varepsilon_n^2/2)b_n.
\]
This verifies the alignment assumption in both cases.
\end{proof}

\subsubsection{Diluted models: Proof of Corollary~\ref{cor:diluted-tensor-search}}

\begin{proof}[Proof of Corollary~\ref{cor:diluted-tensor-search}]
For those $n$ with $r_{\comp}(n;c)>c^2n^{p-1}$, use the polynomial tests of Theorem~\ref{thm:diluted-polynomial-tests}.
For the remaining $n$, the assumption $r=o(r_{\stat}(n;c))$ gives $r=o(c^2n^{p-1})$, and we construct a net for tensor summands.
For $v\sim\nu^{\otimes n}$, the finite second moment gives
\[
    \frac{\|v\|_2^2}{n}\xrightarrow{\mathbb P}1,
    \qquad
    \frac{\max_i v_i^2}{n}\xrightarrow{\mathbb P}0.
\]
The second assertion follows from
$n\Prx[v_1^2>\varepsilon n]\to0$ for each $\varepsilon>0$.
In the expansion of $(\sum_i v_i^2)^p$, the terms with distinct indices
sum to $p!\|X(v)\|_2^2$.
The sum of the terms with a repeated index, divided by $n^p$, is bounded by
\[
    \frac{\binom p2}{n^p}\left(\sum_i v_i^4\right)
        \left(\sum_i v_i^2\right)^{p-2}
    \le\binom p2\frac{\max_i v_i^2}{n}
        \left(\frac{\sum_i v_i^2}{n}\right)^{p-1}
    \xrightarrow{\mathbb P}0.
\]
Consequently,
\[
    \frac{\|X(v)\|_2^2}{n^p}\xrightarrow{\mathbb P}\frac1{p!}.
\]
Set $b_p\colonequals1/(2\sqrt{p!})$ and define the compact set
\[
    \mathcal U_n\colonequals\left\{
        v\in\mathbb R^n:\|v\|_2\le2\sqrt n,
        \ \|X(v)\|_2\ge b_p n^{p/2}
    \right\}.
\]
Then $\nu^{\otimes n}(\mathcal U_n)\to1$.
For $v,w\in\mathcal U_n$, telescoping the full tensor products gives
\[
    \|X(v)-X(w)\|_2
    \le\|v^{\otimes p}-w^{\otimes p}\|_{\mathrm F}
    \le p(2\sqrt n)^{p-1}\|v-w\|_2.
\]
It follows that
\[
    \left\|\frac{X(v)}{\|X(v)\|_2}
        -\frac{X(w)}{\|X(w)\|_2}\right\|_2
    \le\frac{2p(2\sqrt n)^{p-1}}{b_p n^{p/2}}\|v-w\|_2.
\]
Choose a constant $\delta_p>0$ small enough that
$2^p p\delta_p/b_p\le1$.
A maximal $\delta_p\sqrt n$-separated subset of $\mathcal U_n$ is a
$\delta_p\sqrt n$-net of that set.
Comparing the volumes of disjoint balls of radius
$\delta_p\sqrt n/2$ shows that its cardinality is at most
$(1+4/\delta_p)^n$.
The normalized tensor images of these points therefore form a
$1$-net of the normalized summand directions from $\mathcal U_n$.
The net case of Theorem~\ref{thm:scan-detection} applies to $\mathcal S_n=X(\mathcal U_n)$ with
$a_n=(b_p/2)n^{p/2}$ and $M_n\le\exp(C_p n)$.
Its sufficient condition is $r=o(c^2n^{p-1})$, as required.
If $\nu$ has finite support of size $L$, a brute-force search over all nonzero normalized
tensor images of its $L^n$ possible coordinate vectors gives the same
conclusion by the discrete case, restricting the typical set to
$\mathcal U_n$.
\end{proof}

\section*{Acknowledgments}
\addcontentsline{toc}{section}{Acknowledgments}

Thanks to Dan Mikulincer, Cris Moore, and Alex Wein for helpful discussions over the course of this project.

\phantomsection
\addcontentsline{toc}{section}{\refname}
\bibliographystyle{alpha}
\bibliography{main}

\end{document}